\documentclass[10pt,journal]{IEEEtran}
\ifCLASSOPTIONcompsoc
  \usepackage[nocompress]{cite}
\else
  \usepackage{cite}
\fi
\usepackage{amssymb}
\usepackage{amsmath}
\usepackage{dsfont}
\usepackage{amsthm}
\usepackage{graphicx}
\usepackage{mathrsfs,doi,color}
\usepackage{subfigure}

\newtheoremstyle{mythm}{1.5ex plus 1ex minus .2ex}{1.5ex plus 1ex minus .2ex} {\rm}{\parindent}{\it\it}{\rm{:}}{1em}{}
\theoremstyle{mythm}

\renewcommand{\IEEEQED}{\IEEEQEDopen}

\newtheorem{theorem}{Theorem}[section]
\newtheorem{corollary}{Corollary}[section]
\newtheorem{lemma}{Lemma}[section]
\newtheorem{proposition}{Proposition}[section]
\newtheorem{definition}{Definition}[section]
\newtheorem{remark}{Remark}

\newcommand{\Prob}{\mathbb{P}}

\newcommand{\mcN}{\mathcal{N}}

\ifCLASSINFOpdf
\else
\fi
\begin{document}

\title{Convergence Analysis and Strategy Control of Evolutionary Games with Imitation Rule on Toroidal Grid: A Full Version}

\author{Ge~Chen, Yongyuan Yu
  \IEEEcompsocitemizethanks{\IEEEcompsocthanksitem
This material is supported  by the Strategic Priority
Research Program of Chinese Academy of Sciences (XDA27000000), by
the National Natural Science Foundation of China (62203264,72192800, 12288201,
 12071465), by the National Key Research and Development
Program of China (2022YFA1004600), and by Natural Science Fund of Shandong Province
(ZR2022QF061).

  \IEEEcompsocthanksitem
 Ge Chen is with the Key Laboratory of Systems and
  Control, Academy of Mathematics and Systems Science, Chinese Academy of Sciences, Beijing 100190,
  China, {\tt  chenge@amss.ac.cn}.
  Yongyuan Yu is with the School of Mathematics, Shandong University, Jinan 250100, China, {\tt  yyyu@sdu.edu.cn}
}
}

\IEEEtitleabstractindextext{%
  \begin{abstract}
 This paper investigates discrete-time evolutionary games with a general stochastic imitation rule on the toroidal grid, which is a grid network with periodic boundary conditions. The imitation rule has been considered as a fundamental rule to the field of evolutionary game theory,
while the grid is treated as the most basic network and has been widely used in the research of spatial (or networked) evolutionary games.
However, currently the investigation of evolutionary games on grids
mainly uses simulations or approximation methods, while few strict analysis is carried out on one-dimensional grids. This paper proves the convergence of evolutionary prisoner's dilemma, evolutionary snowdrift game, and
evolutionary  stag hunt game with the imitation rule on the two-dimensional grid, for the first time to our best knowledge.
Simulations show that our results may almost reach the critical convergence condition for the evolutionary snowdrift (or hawk-dove, chicken) game.
Also, this paper provides some theoretical results for the strategy control of evolutionary games, and solves the Minimum Agent Consensus Control (MACC) problem under some parameter conditions.
We show that for some evolutionary games (like the evolutionary prisoner's dilemma) on the toroidal grid, one fixed defection node can drive all nodes
almost surely converging to defection, while at least four fixed cooperation nodes are required to lead all nodes almost surely converging to cooperation.
\end{abstract}

}

\maketitle

\IEEEdisplaynontitleabstractindextext

\IEEEpeerreviewmaketitle

\renewcommand{\thesection}{\Roman{section}}
\section{Introduction}
\renewcommand{\thesection}{\arabic{section}}
Cooperation is  one of most common behaviors in human society and nature.
When Charles Darwin was doing the great work on the origin of species, he was puzzled by the phenomenon that animals generally form social groups in which most individuals work for the common good.
He believed that natural selection could encourage altruistic behavior among relatives, thus improving the fertility potential of the ``family".
Despite more than one century of research,
the details of how and why cooperation evolved remain to be solved. In the 125th
anniversary issue of Science, the magazine listed 125 fundamental scientific problems of the 21st century \cite{DK-CN:2005}. Among these problems, the 25 most important were highlighted, one of which was ``How did cooperative behavior evolve?"

Evolutionary game theory has become a major mathematical tool to quantify cooperative behaviors under different circumstances. Traditional game theory concerns one-shot two-player games,  however it does not always reflect real situations because players live in a social network and may game many rounds with same or different opponents. In fact, the network reciprocity has been considered as one of the main mechanisms accounting for the evolution of cooperation \cite{MAN:2006B}, and spatial (or networked) evolutionary games have been developed in recent decades. A pioneering work, made by Nowak and May \cite{MAN-RMM:1992,MAN-RMM:1993}, proposed an evolutionary prisoner's dilemma  model on the grid and demonstrated how cooperators can resist the invasion of defectors by simulations.
Later, evolutionary prisoner's dilemma  with stochastic rules on the grid  were studied, where the players  have a probability
depending on the payoff difference to adopt one of the neighboring strategies \cite{Nakamaru1997The,Nakamaru1998Score,Santos2017Phase}.
There also exist some evolutionary prisoner's dilemma on grids with some additional rules. For example, several literature took into account the memory effect, in which players can update their
strategy by considering previous payoffs \cite{Qin2008Effect,Liu2010Memory}; 
Chiong and Kirley studied the evolution of cooperation when players randomly move
on toroidal grids  \cite{Chiong2012Random}; Mahmoodi and Grigolini considered the effect of the social pressure
to the player's  choice between cooperation and defection  \cite{Korosh2017Evolutionary}. Meanwhile,
there are also some work in other types of evolutionary games on grids, like the evolutionary snowdrift game (or hawk-dove game, chicken game) \cite{1996Spatial}, evolutionary stag hunt game \cite{Dong2019Memory}, and general evolutionary game \cite{EL-CH-MN:05}.

With the development of spatial evolutionary game models, the theoretical analysis has attracts more and more attention. However, most spatial evolutionary game models exhibit very complicated dynamics, and they are, therefore, generally difficult to analyze \cite{doebeli2005models}. One mainstream method is to approximate evolutionary game dynamics by differential equations which generally assume the population is large-enough and well-mixed \cite{DM-CM:15,Morimoto2016Subsidy,Hofbauer2003EVOLUTIONARY,lasry2007mean,wang2014hierarchical,Wang2019Mean,Moon2017Linear,Riehl2018A}. 
Another important direction is using semi-tensor product method which can transform spatial evolutionary game dynamics into linear systems \cite{guo2013algebraic,cheng2014finite,cheng2015modeling,cheng2016decomposed}.
There are also some well analyzed evolutionary games, like asynchronous evolutionary games \cite{Riehl2018A},  decision-making evolutionary games \cite{Altman2010Markov,barreiro2020discrete}, state-based games \cite{LI2020108615}, and potential games \cite{Marden2012State,cheng2014finite}.
In addition, there exists some analysis of spatial evolutionary games based on some special networks, like the cycle \cite{EL-CH-MN:05,Ohtsuki2006Evolutionary} and
complete graph \cite{tan2016analysis,Riehl2017Towards,starnini2011coordination}.

This paper will study the evolutionary game with  imitation rule on toroidal grids.
The imitation rule is relevant to  the behaviors of animals, simple organisms, and human society \cite{2010Human,Berg2015Focus}, and
has been considered as a fundamental rule to the field of evolutionary game theory \cite{Riehl2018A}. Because of the importance, the imitation dynamics has attracted a lot of attention in the theoretic anslysis \cite{Riehl2018A}.
However, the discrete-time imitation dynamics is highly nonlinear in general,  which  makes  the research on its asymptotic behaviors challengeable.
One the other hand, just summarized by  Perc et al. \cite{Perc2013Evolutionary}, grids represent very simple
topologies, and provide a very useful entry
point for studying the consequences of structure on the evolution
of cooperation. Also, there are some practical systems, especially in biology and ecology, in which competition between species can be fully represented by grids. Generally speaking, a grid can be regarded as an even field of all competitive strategies, given the possibility of network reciprocity.
Thus, the grid  has been widely used in the research of spatial evolutionary games \cite{MAN-RMM:1992,MAN-RMM:1993,Nakamaru1997The,Nakamaru1998Score,Santos2017Phase,Qin2008Effect,Liu2010Memory,
Chiong2012Random,Korosh2017Evolutionary,1996Spatial,EL-CH-MN:05}.

The main contribution of this paper can be formulated by the following two aspects:

First, this paper provides some general convergence conditions for the discrete-time evolutionary game with imitation rule on the toroidal grid, for the first time to our best knowledge. The convergence is one of most basic properties in the research of spatial evolutionary games; however, due to the nonlinearity of dynamics, the general convergence conditions of spatial evolutionary games remain to be discovered \cite{Riehl2018A}. Current, convergence analysis is carried out under some assumptions or conditions which are special or usually not easy to verify. For example, Riehl \emph{et al.} proved the convergence of the imitation dynamics when all players are so called opponent-coordinating agents \cite{Riehl2018A};
some literature explored convergence conditions using the well-known semi-tensor product method  \cite{guo2013algebraic,cheng2015modeling}, however these conditions are hard to verify for large-scale networks; another important convergence result was provided  under an assumption that the evolutionary game has a global potential function \cite{Marden2012State}, however whether an evolutionary game has a global potential function is generally unknown.
Currently, the investigation of evolutionary games on grids
mainly uses simulations or approximation methods \cite{MAN-RMM:1992,MAN-RMM:1993,Nakamaru1997The,Nakamaru1998Score,Santos2017Phase,Qin2008Effect,Liu2010Memory,
Chiong2012Random,Korosh2017Evolutionary,1996Spatial}, while few strict analysis is carried out on one-dimensional grids \cite{EL-CH-MN:05,Ohtsuki2006Evolutionary,foxall2016evolutionary}. This paper proves the convergence of evolutionary prisoner's dilemma, evolutionary snowdrift game, and
evolutionary  stag hunt game with imitation rule on the two-dimensional grid, where our convergence conditions depend on system parameters and initial states only.
Simulations show that our results may almost reach the critical convergence condition for the evolutionary snowdrift (or hawk-dove, chicken) game.

Second, this paper provides some theoretical results for the strategy control of evolutionary game
on the toroidal grid, and solves the Minimum Agent Consensus Control (MACC) problem under some parameter conditions. The MACC problem (Problem 1 in \cite{Riehl2017Towards}) is to find the smallest set of fixed strategy players which can drive all players converging to a desired consensus state. This problem is
complex to solve in general. Riehl and Cao solved the MACC problem for imitation dynamics on complete, star and ring networks respectively  \cite{Riehl2017Towards}. The approximate solution of the MACC problem has been explored for the imitation dynamics on tree networks \cite{RR-CM:14,Riehl2017Towards}, while simulations have studied the evolutionary prisoner's dilemma under imitative dynamics on scale-free networks \cite{DZLC2010}. Another theoretic framework  for the strategy control of  evolutionary game is to use semi-tensor product method\cite{cheng2015modeling}, however, this method can be intractable for large populations. Different from previous work, we show that for some evolutionary games (like the evolutionary prisoner's dilemma) on the  grid, one fixed defection node can drive all nodes
almost surely converging to defection, while at least four fixed cooperation nodes are required to lead all nodes almost surely converging to cooperation.

The reminder of this paper is organized as follows.
Section~\ref{Mod_sec} introduces the evolutionary game with stochastic imitation dynamics on a toroidal grid.
Section~\ref{ConResult} contains our almost surely convergence results and relevant proofs, while
Section~\ref{IntCR} provides the interventions to total defection and cooperation.
Some simulations are provided in Section \ref{Simulations}, which is
followed by some concluding remarks in Section \ref{Conclusions}.

\renewcommand{\thesection}{\Roman{section}}
\section{Stochastic Evolutionary Game Model}\label{Mod_sec}
\renewcommand{\thesection}{\arabic{section}}
\subsection{Toroidal Grid}

Let $N\geq 3$ and $M\geq 3$ be two integers.  Assume there are $N\times M$ grid points in $\mathbb{Z}^2$ whose coordinates are $\{(i,j)\}_{i\in\{1,\ldots,N\},j\in\{1,\ldots,M\}}$.
To overcome the effect of the boundary, we consider a grid network with periodic boundary conditions
which means for any node $(i,j)$ we have $(i\pm N,j)=(i,j)=(i,j\pm M)$.
Two nodes $(i_1,j_1)$ and $(i_2,j_2)$ are neighbors (labeling as $(i_1,j_1)\sim(i_2,j_2)$) if and only if
\begin{equation*}
\min_{k\in\mathbb{Z}}\big\{|i_1-i_2+k N|\big\}
+\min_{k\in\mathbb{Z}}\big\{|j_1-j_2+k M|\big\}=1.
\end{equation*}
Thus, all grid nodes have exact four neighbors.
Let $\mathcal{G}_{N,M}$  be the $N\times M$ toroidal grid.

\subsection{Evolutionary game model with stochastic imitation rule}

This paper studies a basic  evolutionary games model on the toroidal grid $\mathcal{G}_{N,M}$.
Assume each player is denoted by
a grid point  $(i,j)$.
Each player $(i,j)$ has a time-varying strategy $S_{i,j}(t)$ which takes
values $C$(cooperation) or $D$(defection). This paper considers that the game between any two players is symmetric, and has a payoff matrix
\begin{table}[hbp]\label{Tab_0}
\centering
\begin{tabular}{| c | c | c |}
\hline
      & C & D \\
\hline
      C & $p_1$ & $p_2$  \\
\hline
      D & $p_3$ & $p_4$ \\
 \hline
\end{tabular}
\end{table}\\
Note that the above matrix corresponds to the prisoner's dilemma if $p_3>p_1>p_4>p_2$,  the snowdrift game (hawk-dove game, chicken game) if $p_3>p_1>p_2>p_4$, and the stag hunt game if $p_1>p_3>p_4>p_2$ \cite{Roca2009Effect}.
 At each time step $t$, which denotes
one generation of the discrete evolutionary time, each node $(i,j)$
in the network plays with all its neighbors and computes its obtained payoff by
\begin{multline}\label{payfun}
P_{i,j}(t):=\\
\left\{
\begin{aligned}
k p_1+(4-k) p_2, ~\mbox{if $(i,j)$ adopts cooperation strategy} \\
\mbox{and has exact $k$ cooperation neighbors at time $t$,}\\
k p_3+(4-k) p_4, ~\mbox{if $(i,j)$ adopts defection strategy~~} \\
\mbox{and has exact $k$ cooperation neighbors at time $t$.}
\end{aligned}\right.
\end{multline}

Our model adopts the imitation dynamics with a general stochastic rule to update the strategy of each node at each time. In detail,  every node $(i,j)$ independently and uniformly selects one of its neighbor node
$\mathcal{R}_{i,j}(t)$, and compares its payoff $P_{i,j}(t)$ with $P(\mathcal{R}_{i,j}(t))$, which is the payoff of
the node $\mathcal{R}_{i,j}(t)$ at time $t$. If $P_{i,j}(t)\geq P(\mathcal{R}_{i,j}(t))$, the node $(i,j)$ does not change its current strategy, i.e.,
$S_{i,j}(t+1)=S_{i,j}(t)$. Otherwise,  with a positive probability the node $(i,j)$ adopts the
strategy of the node $\mathcal{R}_{i,j}(t)$ for the next step. Throughout this paper the imitation probability satisfies the following assumption:

\textbf{(A1)} There exists a positive constant $\delta>0$ such that
\begin{eqnarray*}\label{updaterule_1}
&&\Prob \left\{S_{i,j}(t+1)=S(\mathcal{R}_{i,j}(t))| P_{i,j}(t)<P(\mathcal{R}_{i,j}(t)) \right\}\\
&&~\geq \delta, ~~~~\forall 1\leq i\leq N, 1\leq j\leq M, t\geq 0.
\end{eqnarray*}

\begin{remark}
The assumption (A1) is the process that each node has a positive probability (can be probability one) to imitate its neighbor with a higher gain. In fact, (A1) contains a wide class of stochastic imitation rules, like the proportional imitation rule \cite{Riehl2017Towards},  the simplified Fermi imitation rule
\cite{cheng2015modeling}, and some other imitation rules \cite{FCS-JMP:2005,Santos2017Phase,starnini2011coordination}.
\end{remark}
To simplify exposition we abbreviate the above stochastic evolutionary game as \emph{SEG}.


\renewcommand{\thesection}{\Roman{section}}
\section{Convergence of SEG}\label{ConResult}
\renewcommand{\thesection}{\arabic{section}}
Before stating our convergence results, we need to define the probability
space. Let $\mathcal{N}:=\{(i,j)\}_{1\leq i\leq N,1\leq j\leq M}$ denote the set of all grid points.
For the SEG, we let $\Omega=(\mcN^{N\times M} \times \{C,D\}^{N\times M})^{\infty}$ be the sample space, $\mathcal{F}$ be
the Borel $\sigma$-algebra of $\Omega$, and $\Prob$ be the probability measure
on $\mathcal{F}$. Then the probability space of the SEG is written as
$(\Omega,\mathcal{F},\Prob)$.

Let $S(t):=[S_{i,j}(t)]_{i\in\{1,\ldots,N\},j\in\{1,\ldots,M\}}\in \{C,D\}^{N\times M}$ be the strategy matrix of all nodes at time $t$. We first show that the evolutionary prisoner's dilemma will converge to fixed strategies a.s. in finite time.

\begin{theorem}[Convergence of evolutionary prisoner's dilemma]
  \label{Main_result}
  Consider the SEG on $\mathcal{G}_{N,M}$ satisfying $p_3>p_1>p_4>p_2$ and (A1).
For any  initial  strategies,  there exists a finite time $T$ a.s.  such that $S(t)=S(T)$ for all $t\geq T$. In addition, if
 there exist two adjacent nodes $(i,j)\sim (k,l)$ with $S_{i,j}(T)\neq S_{k,l}(T)$, then
 $P_{i,j}(T)=P_{k,l}(T)$.
\end{theorem}

The SEG is highly nonlinear and hard to analyze.  We adopt
the method of ``transforming the analysis of a stochastic system into the
design of control algorithms" first proposed by \cite{GC:17b}, and also used in \cite{GC-WS-SD-YH:19,GC-WS-WM-FB:20,WS-XC-YY-GC:2019}. This method
requires the construction of a new system called as controllable evolutionary game (CEG) to analyze the SEG.
We will introduce this method in Subsection \ref{sec_method}, and put the proof of Theorem \ref{Main_result} in Subsection \ref{sec_proof_1}.

\begin{remark}
The convergence of the evolutionary prisoner's dilemma on grid networks has attracts lots of attention \cite{MAN-RMM:1992,MAN-RMM:1993,Nakamaru1997The,Nakamaru1998Score,Santos2017Phase}, however the theoretic analysis is challenging. Currently, the theoretical research mainly approximates the dynamics to some differential equations in which the population is assumed to be large-enough and well-mixed \cite{Riehl2017Towards}, while, a few exact analysis is carried out on one-dimensional grids \cite{EL-CH-MN:05,Ohtsuki2006Evolutionary,foxall2016evolutionary}, or under some special assumptions
or conditions \cite{Riehl2018A,guo2013algebraic,cheng2015modeling,Marden2012State}. To our best knowledge, Theorem \ref{Main_result} gives a general and clear convergence condition of the evolutionary prisoner's dilemma on a two-dimensional network for the first time.
\end{remark}

\begin{theorem}[Convergence of evolutionary snowdrift game]
  \label{Main_result2}
  Consider the SEG on $\mathcal{G}_{N,M}$ satisfying $p_3>p_1>p_2>p_4$ and (A1). Assume that $p_1+p_2< p_3+p_4$ and $4p_2<p_3+3 p_4$. Then,
for any  initial  strategies,  there exists a finite time $T$ a.s.  such that $S(t)=S(T)$ for all $t\geq T$. In addition, if
 there exist two adjacent nodes $(i,j)\sim (k,l)$ with $S_{i,j}(T)\neq S_{k,l}(T)$, then
 $P_{i,j}(T)=P_{k,l}(T)$.
\end{theorem}

The proof of Theorem \ref{Main_result2} uses the similar idea in the proof of Theorem \ref{Main_result}, and is postponed to Appendix \ref{sec_proof_mr2}.

\begin{remark}
The traditional two-player snowdrift game has the payoff matrix
\begin{table}[hbp]\label{Tab_a0}
\centering
\begin{tabular}{| c | c | c |}
\hline
      & C & D \\
\hline
      C & $b-c/2$ & $b-c$  \\
\hline
      D & $b$ & $0$ \\
 \hline
\end{tabular}
\end{table}\\
with constants $0<c<b$. It can be verified that if $c<b<4c/3$  then the conditions concerning $p_1,\ldots,p_4$ in Theorem \ref{Main_result2} are satisfied.

The traditional two-player hawk-dove game has the payoff matrix
\begin{table}[hbp]\label{Tab_a0}
\centering
\begin{tabular}{| c | c | c |}
\hline
      & C & D \\
\hline
      C & $b/2$ & $0$  \\
\hline
      D & $b$ & $(b-c)/2$ \\
 \hline
\end{tabular}
\end{table}\\
with constants $0<b<c$. It can be verified that if $b<c<5b/3$  then the conditions concerning $p_1,\ldots,p_4$ in Theorem \ref{Main_result2} are satisfied.

The traditional two-player chicken game has the payoff matrix
\begin{table}[hbp]\label{Tab_a0}
\centering
\begin{tabular}{| c | c | c |}
\hline
      & C & D \\
\hline
      C & $b/2$ & $0$  \\
\hline
      D & $b$ & $-c$ \\
 \hline
\end{tabular}
\end{table}\\
with constants $b>0,c>0$. It can be verified that if $b>3c$  then the conditions concerning $p_1,\ldots,p_4$ in Theorem \ref{Main_result2} are satisfied.

Simulations in Section \ref{Simulations} show that the above relations $b<4c/3$, $c<5b/3$ and $b>3c$ almost reach the critical conditions for the convergence of evolutionary snowdrift, hawk-dove, and chicken games
respectively.
\end{remark}

Let $S_C:=[C]_{N\times M}$ be the strategy matrix whose entries are all equal to $C$,  while
$S_D:=[D]_{N\times M}$ be the strategy matrix whose entries are all equal to $D$. By Theorems \ref{Main_result} and \ref{Main_result2} we can get the following corollary:

\begin{corollary}[Convergence to total cooperation or defection]
  \label{cor_1}
  Consider the SEG on $\mathcal{G}_{N,M}$, and assume that the conditions in Theorem \ref{Main_result}
or Theorem \ref{Main_result2} are satisfied.
 If
\begin{multline}\label{cor1_0}
\{p_1+3p_2,2p_1+2p_2,3p_1+p_2\}\cap \{p_3+3p_4,2p_3+2p_4\}\\
=\emptyset,
 \end{multline}
then for any  initial  strategies there exists a finite time $T$ a.s.  such that $S(T)\in \{S_C,S_D\}$.
\end{corollary}
\begin{IEEEproof}
For any time $t$ and any adjacent nodes $(i,j)\sim (k,l)$, if $S_{i,j}(t)\neq S_{k,l}(t)$,  by (\ref{payfun}), (\ref{cor1_0}), and the conditions in Theorem \ref{Main_result}
or Theorem \ref{Main_result2},  we have $P_{i,j}(t)\neq P_{k,l}(t)$.
By Theorems \ref{Main_result} and \ref{Main_result2}, the convergence state must be total cooperation or total defection.
\end{IEEEproof}

\begin{theorem}[Convergence of evolutionary stag hunt game to total cooperation]
  \label{Main_result3}
  Consider the SEG on $\mathcal{G}_{N,M}$ satisfying $p_1>p_3>p_4>p_2$ and (A1). Assume that $p_1+p_2>2p_3$. If there are four cooperation nodes forming a square at the initial time, i.e., there exists
$(i,j)$ satisfying $S_{i,j}(0)=S_{i+1,j}(0)=S_{i,j+1}(0)=S_{i+1,j+1}(0)=C$, then a.s. the strategies of all nodes will converge to cooperation in finite time.
\end{theorem}

\begin{remark}
Some papers assume the payoff matrix in the  evolutionary stag hunt game   by
\begin{table}[hbp]\label{Tab_a0}
\centering
\begin{tabular}{| c | c | c |}
\hline
      & C & D \\
\hline
      C & $1$ & $-r$  \\
\hline
      D & $r$ & $0$ \\
 \hline
\end{tabular}
\end{table}\\
with constant $r\in(0,1)$ \cite{starnini2011coordination,Dong2019Memory}. It can be verified that if $r<1/3$  then the conditions concerning $p_1,\ldots,p_4$ in Theorem \ref{Main_result3} are satisfied.
Simulations in Section \ref{Simulations} show that  if $r\geq 1/3$ the evolutionary stag hunt game can still converge to total cooperation in finite time, which means Theorem \ref{Main_result3} may be further improved.
\end{remark}

Theorems \ref{Main_result},  \ref{Main_result2}, and \ref{Main_result3} give convergence results of SEG under
some constraints concerning $p_1,\ldots,p_4$. The necessary and sufficient condition of $p_1,\ldots,p_4$ for the convergence of SEG is left to the future.

\textbf{Open Problem:}
  Consider the SEG on $\mathcal{G}_{N,M}$ satisfying (A1). What is the necessary and sufficient condition concerning $p_1,\ldots,p_4$ such that
  SEG converges to a fixed state a.s. in finite time under any initial strategies?

\subsection{CEG and connection to SEG}\label{sec_method}

The CEG is a controllable deterministic version of SEG where the stochastic items in SEG
are replaced by control inputs. In detail,  there are still $N\times M$ players on $\mathcal{G}_{N\times M}$, and at each time step $t$, each node $(i,j)$ still accumulates
the obtained payoff $P_{i,j}(t)$.
After that, every node $(i,j)$  updates its state synchronously by arbitrarily picking up a neighbor
$\mathcal{C}_{i,j}(t)$. In other words, the choice of $\mathcal{C}_{i,j}(t)$ is treated as a control input. The update rule of
 $(i,j)$'s  strategy is to compare its payoff $P_{i,j}(t)$ with $P(\mathcal{C}_{i,j}(t))$, which is the payoff of the node $\mathcal{C}_{i,j}(t)$
 at time $t$. If
$P_{i,j}(t)\geq P(\mathcal{C}_{i,j}(t))$, the node $(i,j)$ does not change its strategy for
the next generation. Otherwise, the node $(i,j)$ adopts the current strategy of the node $\mathcal{C}_{i,j}(t)$
for the next generation. We call such dynamics as \emph{controllable evolutionary game}
which is abbreviated by  \emph{CEG}.

 Recall that $S(t)=[S_{i,j}(t)]\in \{C,D\}^{N\times M}$ is the strategy matrix of all nodes at time $t$. Let $A\subseteq \{C,D\}^{N\times M}$ be a set of strategy matrices. We say $A$ \emph{is reached at time $t$} if
$S(t) \in A$, and $A$  \emph{is reached in the time interval $[t_1,t_2]$} if
there exists $t\in [t_1,t_2]$ such that $S(t)\in A$.

\begin{definition}\label{def_reach}
  Let $A_1, A_2\subseteq \{C,D\}^{N\times M}$ be two sets of strategy matrices.  Under the CEG, $A_1$ is said to be \emph{(uniformly) finite-time reachable from
    $A_2$} if there exists a finite duration $t^*>0$ such that for any $S(0) \in
  A_2$, we can find a sequence of control inputs $[\mathcal{C}_{i,j}(0)]_{1\leq i\leq N,1\leq j\leq M}$, $[\mathcal{C}_{i,j}(1)]_{1\leq i\leq N,1\leq j\leq M}$, $\ldots,$ $[\mathcal{C}_{i,j}(t^*-1)]_{1\leq i\leq N,1\leq j\leq M}$ which guarantees $A_1$ is reached in the time interval $[0,t^*]$.
\end{definition}

Based on these definitions we can get the following result.

\begin{lemma}[Connection between SEG and CEG]\label{robust}
  Let $A_1,A_2\subseteq \{C,D\}^{N\times M}$ be two sets of strategy matrices.
Suppose that $A_2$ is an invariant set under SEG, i.e., if $S(0)\in A_2$ then $S(t)\in A_2$ deterministically for any $t\geq 1$.
Assume that $A_1$ is finite-time reachable from $A_2$ under CEG.
Then, under the SEG satisfying (A1), for any $S(0)\in A_2$, $A_1$ is reached in finite time a.s.
\end{lemma}

\begin{IEEEproof}
For any given $t\geq 0$, since $S(t)\in A_2$, and
$A_1$ is finite-time reachable from $A_2$ under the CEG, there exists a sequence of control inputs
  $[\mathcal{C}_{i,j}(t)]$, $[\mathcal{C}_{i,j}(t+1)]$, $\ldots,$ $[\mathcal{C}_{i,j}(t+t^*-1)]$ such that $A_1$ is reached in $[t,t+t^*]$. Let $S'(t+1),S'(t+2),\ldots,S'(t+t^*)$ be the strategy matrices  under the CEG with  controls  $[\mathcal{C}_{i,j}(t)]$, $\ldots,$ $[\mathcal{C}_{i,j}(t+t^*-1)]$.
  Then, under the SEG,  we obtain
  \begin{equation}\label{robust_1}
    \begin{aligned}
      \Prob&\left(\left\{\mbox{$A_1$ is reached in $[t,t+t^*]$}\right\}|S(t)\right)\\
      &\geq \Prob\left(\bigcap_{s=t+1}^{t+t^*} \big\{S(s)=S'(s)\big\}|S(t)\right)\\
      &= \Prob\left(S(t+1)=S'(t+1)|S(t)\right)  \\
      & \times \prod_{s=t+2}^{t+t^*} \Prob\left(S(s)=S'(s)|S(t), \bigcap_{l=t+1}^{s-1} \big\{S(l)=S'(l)\big\} \right).
    \end{aligned}
  \end{equation}
Set the function $H_{i,j}(s)$ to be $\Prob\{S_{i,j}(s+1)=S(\mathcal{C}_{i,j}(s))\}$ if
$P_{i,j}(s)<P(\mathcal{C}_{i,j}(s))$ and $S_{i,j}(s)\neq S(\mathcal{C}_{i,j}(s))$ (the strategy of the node $\mathcal{C}_{i,j}(s)$ at time $s$), and to be $1$ otherwise.
By (A1) we have
$H_{i,j}(s)\geq \delta.$
 Also, by the definitions of SEG and  CEG  we have
   \begin{eqnarray}\label{robust_2}
  &&\Prob\left(S(s)=S'(s)|S(t), \bigcap_{l=t+1}^{s-1} \big\{S(l)=S'(l)\big\} \right)\nonumber\\
  &&~=\Prob\left(S(s)=S'(s)|S(s-1)=S'(s-1)\big\} \right)\nonumber\\
  &&~\geq \prod_{1\leq i\leq N} \prod_{1\leq j\leq M} \big( \Prob\left\{\mathcal{R}_{i,j}(s-1)=\mathcal{C}_{i,j}(s-1)\right\}\nonumber\\
  &&~~~~ \times H_{i,j}(s-1) \big)\nonumber\\
  &&~\geq \left( \frac{\delta}{4} \right)^{NM}.
  \end{eqnarray}
Substituting (\ref{robust_2}) into (\ref{robust_1}) obtains
  \begin{equation}\label{robust_3}
      \Prob\left(\left\{\mbox{$A_1$ is reached in $[t,t+t^*]$}\right\}|S(t)\right)\geq \left( \frac{\delta}{4} \right)^{NMt^*}.
  \end{equation}
Set $E_t$ to be the event that $A_1$ is reached in $[t,t+t^*]$, and
let $E_{t}^c$ be the complement set of $E_{t}$.
For any integer $K>0$ and initial strategy $S(0)$,
Bayes' Theorem and equation~(\ref{robust_3}) imply
\begin{eqnarray*}\label{robust_4}
&&\Prob\big(\left\{\mbox{$A_1$ is not reached in $[0,(t^*+1)K-1]$}\right\}|S(0)\big)\nonumber\\
&&=\Prob\Big(\bigcap_{m=0}^{K-1}E_{m(t^*+1)}^c\big|S(0)\Big)\nonumber\\
&&=\Prob\left(E_{0}^c|S(0)\right)\prod_{m=1}^{K-1}\Prob\Big(E_{m(t^*+1)}^c\big|S(0),\!\!\bigcap_{0\leq m'<m}\!\!E_{m'(t^*+1)}^c\Big)\nonumber\\
&&\leq \Big(1-\left( \frac{\delta}{4} \right)^{NMt^*} \Big)^K.
\end{eqnarray*}
By the Borel-Cantelli lemma $A_1$ is reached in finite time a.s.
\end{IEEEproof}

Lemma \ref{robust} reveals that, to prove the convergence of the stochastic dynamic model,
it suffices to design suitable control such that a
convergence set is reached.
On the basis of the lemma, we give the proof of Theorem~\ref{Main_result}.

\subsection{Proof of Theorem~\ref{Main_result}}\label{sec_proof_1}
Let $\Omega^* \subseteq \{C,D\}^{N\times M}$ be the set of strategy matrices satisfying that two adjacent nodes with different strategies have a same payoff.
This means, if $S(t)\in\Omega^*$ we have
\begin{equation}\label{mrb_1}
(i,j)\sim (k,l), S_{i,j}(t)\neq S_{k,l}(t) \Rightarrow P_{i,j}(t)=P_{k,l}(t).
\end{equation}
By (\ref{mrb_1}) we can get
\begin{equation}\label{mrb_2}
S(t)\in\Omega^* \Rightarrow S(t)=S(t+1)=S(t+2)=\cdots
\end{equation}
under SEG. Thus, if $\Omega^*$ is finite-time
  reachable from any initial strategies under CEG,  by Lemma \ref{robust} and
  (\ref{mrb_2}), our result is obtained.

Therefore, it suffices to design a control algorithm such that $\Omega^*$ is finite-time
  reachable from any given initial $S(0)\notin\Omega^*$ under the CEG.
The algorithm is established by the following three steps:

\textbf{Step 1}: For $t=0,1,\ldots,$ and each node $(i,j)$, the control input $\mathcal{C}_{i,j}(t)$
is selected as follows:
\begin{equation}\label{mresult_1}
\mathcal{C}_{i,j}(t)=
\left\{
\begin{aligned}
\mathop{\arg\max}\limits_{(k,l)\sim(i,j),S_{k,l}(t)=D}P_{k,l}(t), ~~~~\mbox{if $(i,j)$ has} \\
~~~~~~~~~~ \mbox{at least one defection neighbor at $t$,}\\
\mathop{\arg\min}\limits_{(k,l)\sim(i,j)}P_{k,l}(t),~~~~~~~~~~~~~~~~ \mbox{otherwise}.
\end{aligned}\right.
\end{equation}
Here we set $\mathcal{C}_{i,j}(t)$ to be arbitrary one neighbor of $(i,j)$ with maximal (or minimal) payoff
if the mapping $\arg\max$ (or $\arg\min$) in (\ref{mresult_1}) is not unique.

If the node $(i,j)$ is an isolated defection node at time $t$, i.e.,
 $S_{i,j}(t)=D$, and $S_{k,l}(t)=C$ for all $(k,l)\sim (i,j)$, which indicates the node $(i,j)$ has the largest payoff $4p_3$ by (\ref{payfun}) and the condition $p_3>p_1>p_4>p_2$, then under the  CEG and the control input (\ref{mresult_1}) we have
 $S_{i,j}(t+1)=D$, and $S_{k,l}(t+1)=D$ for all $(k,l)\sim (i,j)$.
Then, by (\ref{mresult_1}), for each node $(i,j)$ it is clear that
\begin{equation}\label{mresult_1_2}
S_{i,j}(t)=D  \Rightarrow S_{i,j}(t+1)=D,
\end{equation}
and
\begin{multline}\label{mresult_1_3}
S_{i,j}(t)=C, P_{i,j}(t) < \max_{(k,l)\sim(i,j),S_{k,l}(t)=D}P_{k,l}(t)\\
 \Rightarrow S_{i,j}(t+1)=D.
\end{multline}
If a cooperation node $(i,j)$ has at most one cooperation neighbor at time $t$, then by
(\ref{payfun}), (\ref{mresult_1_3}),  and the relation  $p_3>p_1>p_4>p_2$
  we have
\begin{multline}\label{mresult_1_4}
P_{i,j}(t)\leq p_1+3 p_2 < p_3+3 p_4
 \Rightarrow S_{i,j}(t+1)=D.
\end{multline}

We repeatedly carry out the CEG with the control input (\ref{mresult_1}) until the strategies of all nodes keep unchanged. We record the stop time as $T_1$. Let $n_C(t)$ be the number of cooperation nodes at time $t$.
By (\ref{mresult_1_2})-(\ref{mresult_1_3}) it can be found
\begin{equation}\label{mresult_2}
n_C(0)>n_C(1)>\cdots>n_C(T_1).
 \end{equation}
If $n_C(T_1)=0$, then our result is obtained. Otherwise,
it can be deduced from (\ref{mresult_1_3}) that
 \begin{multline}\label{mresult_1_3b}
(i,j)\sim (k,l), S_{i,j}(T_1)=C,  S_{k,l}(T_1)=D\\
 \Rightarrow  P_{i,j}(T_1) \geq P_{k,l}(T_1).
\end{multline}
Also, by (\ref{mresult_1_4}), each cooperation node has at least two cooperation neighbors at time $T_1$.
Next we carry out the following Step 2.

\textbf{Step 2}:
If there exist some cooperation nodes which have exact $2$ cooperation neighbors at time $T_1$, and have different payoff from defection neighbors,
 without loss of generality we assume $(i,j)$ has exact $2$ cooperation neighbors, and
\begin{equation}\label{mrb_3}
\left\{
\begin{aligned}
&S_{i,j}(T_1)=S_{i+1,j}(T_1)=C, S_{i,j+1}(T_1)=D, \\
&P_{i,j}(T_1)\neq  P_{i,j+1}(T_1).
\end{aligned}\right.
\end{equation}
By (\ref{payfun}), (\ref{mrb_3}) and (\ref{mresult_1_3b}) we have
\begin{equation}\label{mrb_3_1}
P_{i,j}(T_1)=2p_1+2p_2 > P_{i,j+1}(T_1).
\end{equation}
By (\ref{payfun}), (\ref{mrb_3_1}), and the relation $p_3>p_1>p_4>p_2$ we get $P_{i,j+1}(T_1)=p_3+3p_4<2p_1+2p_2$ (If $p_3+3p_4\geq 2p_1+2p_2$ the Step 2 can be skipped). Thus, the node $(i,j+1)$ has only one cooperation neighbor at time $T_1$, which indicates  $S_{i+1,j+1}(T_1)=D$. Also, by (\ref{mresult_1_3}) and (\ref{payfun}) again we have
\begin{equation}\label{mrb_3_2}
P_{i+1,j+1}(T_1)\leq P_{i+1,j}(T_1)\leq 3p_1+p_2.
\end{equation}
Next we  discuss (\ref{mrb_3_2}) with two cases.

Case I: $P_{i+1,j+1}(T_1)=P_{i+1,j}(T_1)$. Because $p_3+3p_4<2p_1+2p_2<2p_3+2p_4$,  we have  $P_{i+1,j+1}(T_1)=3p_1+p_2$.  Choose $\mathcal{C}_{i,j+1}(T_1)=(i,j)$ for node $(i,j+1)$, while for any other node we choose a neighbor which has the same strategy with it as the control input at $T_1$ (We note that at time $T_1$ there is no isolated cooperation or detection node by the discussion of Step 1).  Then, by the CEG, at time $T_1+1$ the strategy of the node $(i,j+1)$ changes from $D$ to $C$, while the other nodes keep unvaried. Thus,
we have
\begin{multline}\label{mrb_7}
P_{i+1,j+1}(T_1+1)>P_{i+1,j+1}(T_1)=3p_1+p_2\\
\geq \max\{P_{i,j+1}(T_1+1),P_{i+1,j}(T_1+1)\}.
\end{multline}
We adopt the control input (\ref{mresult_1}) at time $T_1+1$. By (\ref{mrb_7})
and (\ref{mresult_1_3}) we can get
\begin{equation}\label{mrb_8}
S_{i+1,j+1}(T_1+2)=S_{i,j+1}(T_1+2)=S_{i+1,j}(T_1+2)=D.
\end{equation}
 \begin{figure}
  \centering
  \includegraphics[width=3in]{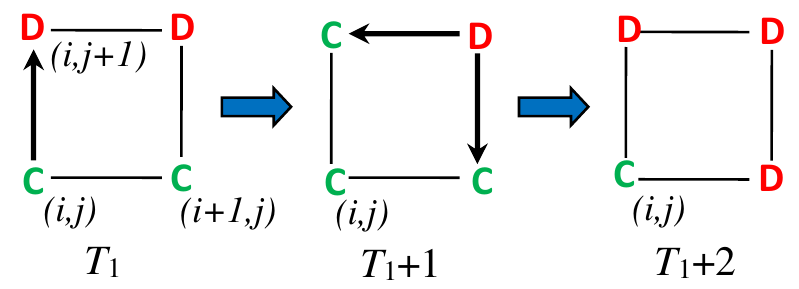}
  \caption{The evolution of nodes' strategies from (\ref{mrb_3}) to (\ref{mrb_8}).}\label{LFig3}
\end{figure}
The evolution of above process is shown in Fig. \ref{LFig3}.
Combining (\ref{mrb_8}) with (\ref{mresult_1_2}) and (\ref{mrb_3}) yields
\begin{equation}\label{mrb_6}
n_{C}(T_1)> n_C(T_1+2).
\end{equation}

Case II: $P_{i+1,j+1}(T_1)<P_{i+1,j}(T_1)$. The defection node $(i+1,j+1)$ must have one or two cooperation neighbor nodes at time $T_1$.
If $(i+1,j+1)$ has two cooperation neighbor nodes at time $T_1$, we choose the same control strategy as Case I. Then
at time $T_1+1$ the strategy of the node $(i,j+1)$ changes from $D$ to $C$, while the other nodes keep unvaried.
So, the defection node $(i+1,j+1)$ has three cooperation neighbor nodes at time $T_1+1$,
which indicates
\begin{multline}\label{mrb_a1}
P_{i+1,j+1}(T_1+1)=3p_3+p_4>3p_1+p_2\\
\geq \max\{P_{i,j+1}(T_1+1),P_{i+1,j}(T_1+1)\}
\end{multline}
by (\ref{payfun}) and the relation $p_3>p_1>p_4>p_2$. We adopt the control input (\ref{mresult_1}) at time $T_1+1$. By (\ref{mrb_a1})
and (\ref{mresult_1_3}) we can get (\ref{mrb_8}), and then (\ref{mrb_6}) still holds.

\begin{figure}
  \centering
  \includegraphics[width=3in]{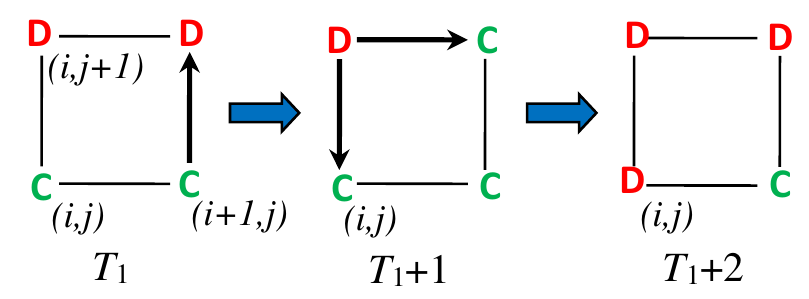}
  \caption{The evolution of nodes' strategies for the case when $(i+1,j+1)$ has one cooperation neighbor at time $T_1$.}\label{LFig2}
\end{figure}

If $(i+1,j+1)$ has one cooperation neighbor node at time $T_1$, we choose $\mathcal{C}_{i+1,j+1}(T_1)=(i+1,j)$ for node $(i+1,j+1)$, while for any other node we choose a neighbor which has the same strategy with it as the control input.  Then, under the CEG, by (\ref{mrb_3}) and (\ref{mrb_3_1}), at time $T_1+1$ the strategy of the node $(i+1,j+1)$ changes from $D$ to $C$, while the other nodes keep strategy unaltered. At time $T_1+1$, the cooperation node $(i+1,j+1)$ still has one cooperation neighbor, while the defection node $(i,j+1)$ has two cooperation neighbors. Thus,
 by (\ref{mrb_3}), (\ref{mrb_3_1}), (\ref{payfun}) and the relation $p_3>p_1>p_4>p_2$ we have
\begin{eqnarray}\label{mrb_a2}
\begin{aligned}
&P_{i,j+1}(T_1+1)=2p_3+2p_4\\
&~~>2p_1+2p_2=P_{i,j}(T_1+1)\\
&~~>p_1+3 p_2=P_{i+1,j+1}(T_1+1).
\end{aligned}
\end{eqnarray}
We adopt the control input (\ref{mresult_1}) at time $T_1+1$. According to (\ref{mrb_a2})
and (\ref{mresult_1_3}) we can get
\begin{equation*}\label{mrb_a3}
S_{i,j}(T_1+2)=S_{i,j+1}(T_1+2)=S_{i+1,j+1}(T_1+2)=D,
\end{equation*}
and then (\ref{mrb_6}) still holds.
Fig. \ref{LFig2} shows the evolution of nodes' strategies for the above case that $(i+1,j+1)$ has one cooperation neighbor at time $T_1$.

For $t\geq T_1+2$  we repeatedly carry out Step 1 and the above process until  the strategies of all nodes keep unchanged under the control input (\ref{mresult_1}), and
 there is no cooperation node
which has $2$ cooperation neighbors and has a payoff different from defection neighbors.
We record the stop time as $T_2$.
If $S(T_2)\in\Omega^*$, our result is obtained. Otherwise, we continue the following Step 3.

\textbf{Step 3}:
Because $S(T_2)\notin\Omega^*$, there are some cooperation nodes which have $3$ cooperation neighbors and have a payoff different from defection neighbors. Without loss of generality we assume
 $S_{i,j}(T_2)=C$, $S_{i,j+1}(T_2)=D$, and
\begin{equation}\label{mrb_9}
P_{i,j}(T_2)=3p_1+p_2 \neq  P_{i,j+1}(T_2).
\end{equation}
By above assumption we get
\begin{equation}\label{mrb_9_1}
S_{i-1,j}(T_2)=S_{i+1,j}(T_2)=C.
\end{equation}
Also, by (\ref{payfun}), (\ref{mrb_9}) and (\ref{mresult_1_3b}) we have
\begin{equation}\label{mrb_9_2}
P_{i,j}(T_2)=3p_1+p_2>P_{i,j+1}(T_2).
\end{equation}
 We continue our discussion with two cases as follows.

\begin{figure}
  \centering
  \includegraphics[width=3in]{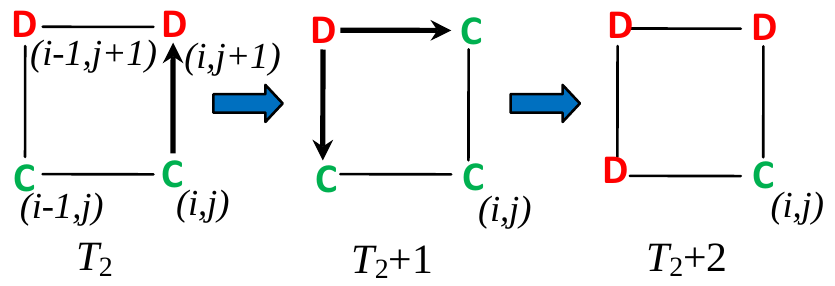}
  \caption{The evolution of nodes' strategies from (\ref{mrb_10}) to (\ref{mrb_15}).}\label{LFig4}
\end{figure}

Case I: $S_{i-1,j+1}(T_2)=D$ and $P_{i-1,j+1}(T_2)=P_{i-1,j}(T_2)$, or, $S_{i+1,j+1}(T_2)=D$ and $P_{i+1,j+1}(T_2)=P_{i+1,j}(T_2)$. Without loss of generality we assume
\begin{equation}\label{mrb_10}
S_{i-1,j+1}(T_2)=D,~~P_{i-1,j+1}(T_2)=P_{i-1,j}(T_2).
\end{equation}
For the node $(i,j+1)$ we choose
$\mathcal{C}_{i,j+1}(T_2)=(i,j)$, while for any other node we choose a neighbor which has the same strategy with it as the control input of $T_2$.
By (\ref{mrb_9_2}),   the strategy of the node $(i,j+1)$ changes from $D$ to $C$, while the other nodes keep strategy unaltered at time $T_2+1$ under the CEG. Thus, by (\ref{mrb_9_2}) and (\ref{mrb_10}),
\begin{equation}\label{mrb_10_a}
S_{i-1,j+1}(T_2+1)=D,~S_{i-1,j}(T_2+1)=S_{i,j+1}(T_2+1)=C.
\end{equation}
Also, by (\ref{mrb_9_2}), (\ref{payfun}) and the relation $p_3>p_1>p_4>p_2$ we obtain that
 $(i,j+1)$ has at most two cooperation neighbors at time $T_2$, so
 \begin{equation}\label{mrb_10_1}
P_{i,j+1}(T_2+1)\leq 2 p_1+2p_2.
\end{equation}
Meanwhile, the defection node $(i-1,j+1)$ has at least two cooperation neighbors $(i,j+1)$ and $(i-1,j)$ at time $T_2+1$, by  (\ref{payfun}), (\ref{mrb_10_1}) and the relation $p_3>p_1>p_4>p_2$ we have
 \begin{eqnarray}\label{mrb_10_2}
 \begin{aligned}
P_{i-1,j+1}(T_2+1)&\geq 2 p_3+2p_4\\
&>2 p_1+2p_2\geq P_{i,j+1}(T_2+1).
\end{aligned}
\end{eqnarray}
On the other hand, because $(i,j+1)$ changes from $D$ to $C$ at time $T_2+1$, by (\ref{mrb_10}) we have
 \begin{eqnarray}\label{mrb_10_3}
 \begin{aligned}
P_{i-1,j+1}(T_2+1)&=P_{i-1,j+1}(T_2)+p_3-p_4\\
&>P_{i-1,j+1}(T_2)=P_{i-1,j}(T_2).
\end{aligned}
\end{eqnarray}
We adopt the control input (\ref{mresult_1}) at time $T_2+1$. On the basis of (\ref{mrb_10_a}), (\ref{mrb_10_2}) and (\ref{mrb_10_3}), it follows
\begin{equation}\label{mrb_15}
S_{i-1,j+1}(T_2+2)=S_{i,j+1}(T_2+2)=S_{i-1,j}(T_2+2)=D.
\end{equation}
The evolution of nodes' strategies from (\ref{mrb_10}) to (\ref{mrb_15}) is shown in Fig. \ref{LFig4}.
Combining (\ref{mresult_1_2}), (\ref{mrb_10}), and (\ref{mrb_15}) yields
\begin{equation}\label{mrb_13}
n_{C}(T_2)> n_C(T_2+2).
\end{equation}

\begin{figure*}
  \centering
  \includegraphics[width=6.5in]{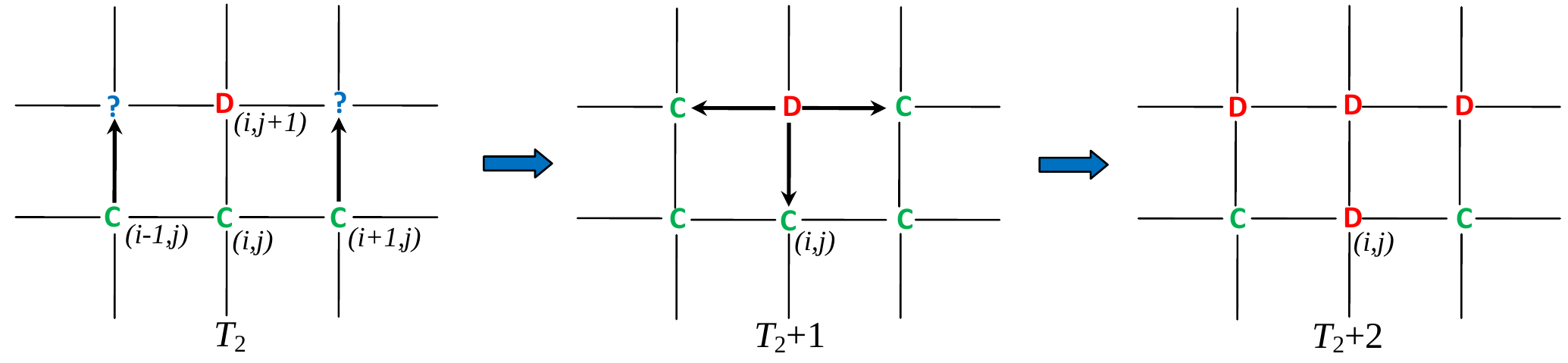}
  \caption{The evolution of nodes' strategies for the Case II of the Step 3 in the proof of Theorem \ref{Main_result}, where  ``?" denotes the unknown strategy.}\label{LFig1}
\end{figure*}

Case II: $S_{i-1,j+1}(T_2)=C$ or $P_{i-1,j+1}(T_2)<P_{i-1,j}(T_2)$, and, $S_{i+1,j+1}(T_2)=C$ or $P_{i+1,j+1}(T_2)<P_{i+1,j}(T_2)$.
For nodes $(i-1,j+1)$ and $(i+1,j+1)$ we choose control inputs
$\mathcal{C}_{i-1,j+1}(T_2)=(i-1,j)$ and $\mathcal{C}_{i+1,j+1}(T_2)=(i+1,j)$, respectively, while for any other node we choose a neighbor which has the same strategy with it as the control input of $T_2$.
Then, under the CEG,  we get
\begin{equation}\label{mrb_14}
S_{i-1,j+1}(T_2+1)=C=S_{i+1,j+1}(T_2+1).
\end{equation}
Meanwhile, the strategies of all nodes except $(i-1,j+1)$ and $(i+1,j+1)$ keep unchanged at time $T_2+1$.
Thus, by the fact that  $S_{i,j}(T_2)=C$, $S_{i,j+1}(T_2)=D$, and $p_3>p_1>p_4>p_2$,  we have
\begin{multline}\label{mresult_2_1b}
P_{i,j+1}(T_2+1)\geq 3p_3+p_4>3p_1+p_2\\
 \geq \max\big\{P_{i,j}(T_2+1),P_{i-1,j+1}(T_2+1),P_{i+1,j+1}(T_2+1)\big\}.
\end{multline}
We adopt the control input (\ref{mresult_1}) at time $T_2+1$.
By resorting to (\ref{mresult_2_1b}), (\ref{mresult_1_2}) and (\ref{mresult_1_3}), we can see
\begin{multline}\label{mresult_2_1c}
S_{i,j+1}(T_2+2)=S_{i,j}(T_2+2)=S_{i-1,j+1}(T_2+2)\\
=S_{i+1,j+1}(T_2+2)=D.
\end{multline}
The evolution of nodes' strategies for Case II is shown in Fig. \ref{LFig1}.
Combining (\ref{mresult_1_2}) and (\ref{mresult_2_1c}) yields (\ref{mrb_13}).

For $t\geq T_2+2$  we repeatedly carry out Steps 1-2 and the above process until $\Omega^*$ is reached. Let $T_3$ be the stop time. By (\ref{mresult_2}), (\ref{mrb_6}), and (\ref{mrb_13}) we have $T_3\leq 2(NM-1)$.

\subsection{Proof of Theorem~\ref{Main_result3}}\label{sec_proof_3}
Without loss of generality we assume $S_{i,j}(0)=S_{i+1,j}(0)=S_{i,j+1}(0)=S_{i+1,j+1}(0)=C$.
Set $\mathcal{I}_0:=\{(i,j),(i+1,j),(i,j+1),(i+1,j+1)\}$.
Let $A_2\subseteq \{C,D\}^{N\times M}$ be the set of strategy matrices satisfying that all nodes in $\mathcal{I}_0$ have cooperation strategy. Thus, we have $S(0)\in A_2$.
By (\ref{payfun}) and the relations $p_1>p_3>p_4>p_2$ and $2p_1+2p_2> 4p_3$,
\begin{multline}\label{SH_1}
\min\left\{P_{i,j}(0),P_{i+1,j}(0),P_{i,j+1}(0),P_{i+1,j+1}(0)\right\}\\
 \geq 2p_1+2p_2>k p_3+(4-k)p_4, ~~0\leq k\leq 4,
\end{multline}
which indicates any defection node have a payoff less than $P_{i,j}(0),P_{i+1,j}(0),P_{i,j+1}(0),P_{i+1,j+1}(0)$. Thus, under SEG,
$S(1)\in A_2$ deterministically. With the similar way we can get $S(t)\in A_2$ deterministically for
$t\geq 2$. By Lemma \ref{robust}, our result is obtained if we can show that $\{S_C\}$ is finite-time reachable from $A_2$ under CEG.

Next we design control inputs under CEG such that $\{S_C\}$ is reached from $A_2$ in finite-time.
Let $\mathcal{I}_1$ be the set of all neighboring nodes of $\mathcal{I}_0$.
At the initial time, for any node $(k,l)\in \mathcal{I}_1\setminus \mathcal{I}_0$
we choose its control input $\mathcal{C}_{k,l}(0)=(k',l')$ with $(k',l')\in \mathcal{I}_0$ and
$(k',l')\sim (k,l)$, while for other nodes we choose arbitrary control inputs.
By (\ref{SH_1}), all nodes in $\mathcal{I}_1$ have cooperation strategy at time $t=1$. Also,
by the definition of $\mathcal{I}_1$, each node in  $\mathcal{I}_1$ has at least two cooperation neighbors. For $s\geq 2$, we set $\mathcal{I}_s$ to be the set of all neighboring nodes of $\mathcal{I}_{s-1}$, and choose the control input $\mathcal{C}_{k,l}(s-1)=(k',l')$ for $(k,l)\in \mathcal{I}_1\setminus \mathcal{I}_{s-1}$, $(k',l')\in \mathcal{I}_{s-1}$,  and
$(k',l')\sim (k,l)$, while for other nodes we choose arbitrary control inputs.
 With the similar discussion as $\mathcal{I}_1$ we can get all nodes in $\mathcal{I}_s$
have cooperation strategy. Because the toroidal grid $\mathcal{G}_{N,M}$ is finite,
there exists a finite $s$ such that $\mathcal{I}_s$ contains all nodes in $\mathcal{G}_{N,M}$.

\renewcommand{\thesection}{\Roman{section}}
\section{Control strategies  of SEG}\label{IntCR}
\renewcommand{\thesection}{\arabic{section}}

An interesting attempt is to intervene the convergence results of spatial evolutionary games.
Currently there exist multiple control methods to affect the dynamics of evolutionary games, in which a most natural and intuitive method
is to fix the strategies of some nodes in the network with the hope of producing a desired global outcome.
Riehl and Cao proposed a Minimum Agent Consensus Control (MACC) problem which is to find the smallest set of constant strategy players which can drive all players converging to a desired consensus state (Problem 1 in \cite{Riehl2017Towards}). Currently, this problem has been  solved exactly on complete, star and ring networks respectively  \cite{Riehl2017Towards}, and solved approximately on tree networks \cite{RR-CM:14,Riehl2017Towards}. In this paper we try to solve the MACC problem on the toroidal grid.

Recall that $\mathcal{N}$ is the node set of $\mathcal{G}_{N,M}$.
Let $\mathcal{N}_C \subset \mathcal{N}$  be the grid nodes
whose strategies are fixed to cooperation, and call the repeated game on $\mathcal{G}_{N,M}$ by \emph{$\mathcal{N}_C$-SEG} if the nodes in $\mathcal{N}\backslash \mathcal{N}_C$
carry out the SEG. Similarly, we set
 $\mathcal{N}_D\subset \mathcal{N}$ to be the grid nodes whose strategies are fixed to defection,
and call the system by
 \emph{$\mathcal{N}_D$-SEG} if the nodes in $\mathcal{N}\backslash \mathcal{N}_D$ implement the SEG.

Similar to the CEG, we define the \emph{$\mathcal{N}_C$-CEG}
(or \emph{$\mathcal{N}_D$-CEG}) if the nodes in $\mathcal{N}\backslash \mathcal{N}_C$
(or $\mathcal{N}\backslash \mathcal{N}_D$) update their strategies by the same way as the
CEG. We first study the MACC problem where the  desired consensus state is defection.


\begin{proposition}[Control to total defection]
\label{Main_result_2}
For any non-empty set  of constant defection nodes $\mathcal{N}_D \subset \mathcal{N}$,
and any initial strategies of $\mathcal{N}\backslash \mathcal{N}_D$, if the parameters $p_1,p_2,p_3,p_4$ satisfy the conditions in Corollary \ref{cor_1}, then the $\mathcal{N}_D$-SEG  converges to total defection in finite time a.s., i.e.,
there exists a finite time $T$ a.s.  such that $S(T)=S_D$.
\end{proposition}
\begin{IEEEproof}
By the similar method to the proofs of Theorems \ref{Main_result} and \ref{Main_result2},
we can get that $\mathcal{N}_D$-SEG  converges to a fixed state in finite time  a.s. By the similar discussion to the proof of Corollary \ref{cor_1} we can get $\mathcal{N}_D$-SEG  converges to total defection in finite time a.s.
\end{IEEEproof}

\begin{remark}
Proposition \ref{Main_result_2} demonstrates that only one constant defection node can cause all nodes to defect with each other eventually. Thus, for the SEG on the toroidal grid, providing the conditions in Corollary \ref{cor_1} concerning the parameters
$p_1,\ldots,p_4$ are satisfied, one node is the solution of the MACC problem with defection being the desired
consensus state.
\end{remark}

Compared with defection being the desired consensus state, the MACC problem with cooperation being the desired consensus state is much complex. We first give some conditions when constant cooperation nodes cannot drive all nodes  converging to cooperation.

\begin{proposition}\label{Main_result_4}
Assume that the number of constant cooperation nodes is less than $4$ in $\mathcal{G}_{N,M}$ with $N\geq 4, M\geq 4$.
If the initial strategy of each node in $\mathcal{N}\backslash \mathcal{N}_C$ is defection, the $\mathcal{N}_C$-SEG with $p_3>p_1>p_4>p_2$  cannot converge to total cooperation  deterministically.
\end{proposition}
\begin{IEEEproof}
Because $|\mathcal{N}_C|\leq 3$,  $N\geq 4, M\geq 4$, and the initial strategies of the nodes in $\mathcal{N}\backslash \mathcal{N}_C$ are defection,  there is at most one cooperation node
which has two cooperation neighbors at the initial time.

If there is no cooperation node which has two cooperation neighbors at the initial time, by $p_3>p_1>p_4>p_2$ and (\ref{payfun}) we have $S(t)=S(0)$  deterministically for all $t\geq 1$.

If there is one cooperation node $(i,j)$ which has two cooperation neighbors at the initial time, by $p_3>p_1>p_4>p_2$ and (\ref{payfun}) it can be computed that $S_{k,l}(t)=D$  deterministically for any $t\geq 1$ and any node $(k,l)$ which is not a neighbor of $(i,j)$.
\end{IEEEproof}

Next we give a sufficient condition, which provides a simple structure of constant cooperation nodes guiding all nodes to cooperation.

\begin{theorem}[Control to total cooperation]
\label{Main_result_3}
Assume that the constant cooperation nodes $\mathcal{N}_C$ form a $N_1\times M_1$ rectangular area in $\mathcal{G}_{N,M}$, i.e.,
$$\mathcal{N}_C=\{(i,j):i_0\leq i<i_0+N_1, j_0\leq j< j_0+M_1\},$$
where $N\geq 5, M\geq 5$ and $N_1\in [2,N-3]\cup\{N\}$ and $M_1\in [2,M-3]\cup\{M\}$ are two integers.
Consider the $\mathcal{N}_C$-SEG  satisfying (A1), (\ref{cor1_0}), and one of the following two conditions:\\
i) $p_3>p_1>p_4>p_2$ and $2p_1+2p_2>p_3+3p_4$;\\
ii) $p_3>p_1>p_2>p_4$ and $2p_3+2p_4>2p_1+2p_2>p_3+3p_4>4p_2$.\\
Then, for any initial strategies of $\mathcal{N}\backslash \mathcal{N}_C$,
the $\mathcal{N}_C$-SEG converges to total cooperation in finite time a.s., i.e.,
there exists a finite time $T$ a.s.  such that $S(T)=S_C$.
\end{theorem}

The proof of Theorem \ref{Main_result_3} is put in Subsection \ref{proof_mr3}.

\begin{remark}
From Proposition \ref{Main_result_4} and Theorem \ref{Main_result_3}, for the evolutionary prisoner's dilemma
on the toroidal grid with $2p_1+2p_2>p_3+3p_4$,  the solution of the MACC problem with cooperation being the desired consensus state is the rectangular area with four nodes. Interestingly,
it is not always true that more constant  cooperation nodes can more easily lead all nodes to cooperation. A simple counter-example is that for the evolutionary prisoner's dilemma on the toroidal grid,  when $|\mathcal{N}_C|=NM-1$, the remaining defection node
will keep defection forever.
\end{remark}

\subsection{Proof of Theorem \ref{Main_result_3}}\label{proof_mr3}
 \begin{figure*}
  \centering
  \includegraphics[width=6.5in]{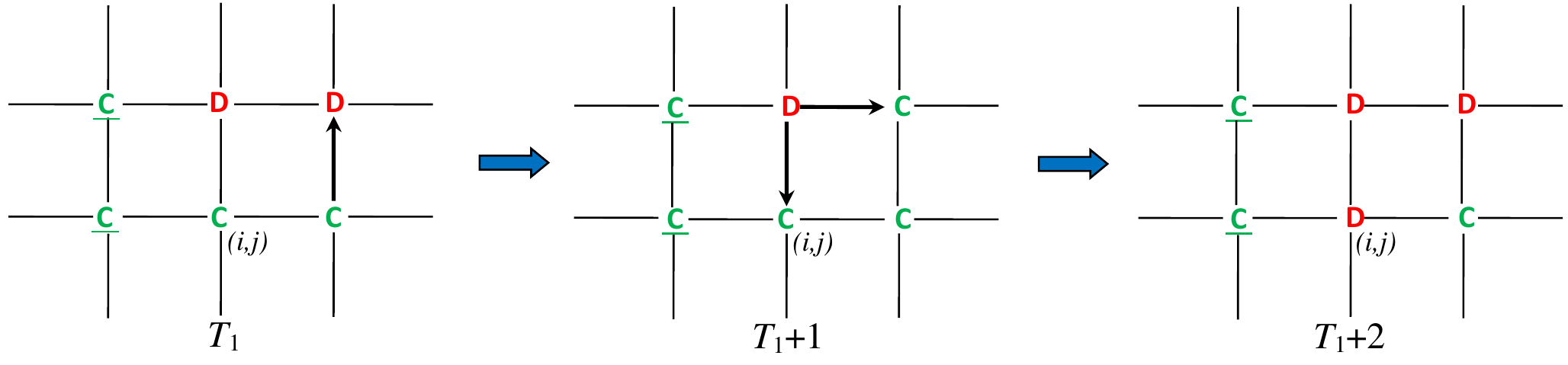}
  \caption{An evolution of nodes' strategies from (\ref{mr3_3}) to (\ref{mr3_5_2}), where ``$\underline{C}$" denotes the node (possibly) belongs to  $\mathcal{N}_C$.}\label{LFig5}
\end{figure*}
By Lemma \ref{robust}, our result is obtained if $\{S_C\}$ is finite-time
  reachable from any initial strategies of $\mathcal{N}\backslash\mathcal{N}_C$ under the $\mathcal{N}_C$-CEG. Thus,
we only need consider the $\mathcal{N}_C$-CEG with $S(0)\neq S_C$, which implies $N_1 M_1<NM$.
Let $S^*$ be the strategy matrix satisfying
\begin{equation}\label{mr3_1}
S_{i,j}^*=
\left\{
\begin{aligned}
C, ~~&\mbox{if~} (i,j)\in \mathcal{N}_C,\\
D, ~~&\mbox{otherwise}.
\end{aligned}\right.
\end{equation}
We first assert that under the $\mathcal{N}_C$-CEG, $\{S^*\}$ is finite-time reachable from
 any initial strategies of $\mathcal{N}\backslash \mathcal{N}_C$.
We carry out the control algorithm by the following $5$ steps:

\textbf{Step 1}: For the nodes in $\mathcal{N}\backslash \mathcal{N}_C$,
we choose (\ref{mresult_1}) as control inputs  until their strategies keep unchanged. We record the stop time as $T_1$. By the similar discussion as the Step 1 of the proof of Theorem~\ref{Main_result} i),  we still obtain (\ref{mresult_2}), and there is no isolated cooperation or defection node at time $T_1$.
If $S(T_1)=S^*$, our assertion holds. Otherwise,
for any nodes $(i,j)\sim (k,l)$ with $(i,j)\notin \mathcal{N}_C$, $S_{i,j}(T_1)=C,$ and $S_{k,l}(T_1)=D$,
similar to (\ref{mresult_1_3b}) we get
\begin{equation}\label{mr3_2}
P_{i,j}(T_1) > P_{k,l}(T_1).
\end{equation}
We carry out the following Step 2.

\textbf{Step 2}:
If there is a cooperation node $(i,j)\notin \mathcal{N}_C$  which has three cooperation neighbors at $T_1$, without loss of generality we assume $S_{i,j+1}(T_1)=D$, and then
\begin{equation}\label{mr3_3}
S_{i-1,j}(T_1)=S_{i,j-1}(T_1)=S_{i+1,j}(T_1)=C.
\end{equation}
Since the nodes in $\mathcal{N}_C$ form a $N_1\times M_1$ rectangular area with $N_1\in [2,N-3]\cup\{N\}, M_1\in [2,M-3]\cup\{M\}$,
the set $\{(i-1,j), (i,j-1), (i+1,j)\}$ has at most one node belongs to $\mathcal{N}_C$.
We first consider the case when $(i-1,j)\in \mathcal{N}_C$. In this case we discuss the strategy of
$(i-1,j+1)$ at time $T_1$.

If $S_{i-1,j+1}(T_1)=C$, because $3p_1+p_2<3 p_3+p_4$  by the conditions of $p_1,\ldots,p_4$,
 using (\ref{mr3_2}) we have
\begin{equation}\label{mr3_4}
3p_1+p_2= P_{i,j}(T_1) > P_{i,j+1}(T_1)=2 p_3+2 p_4.
\end{equation}
By (\ref{mr3_4}) we get $S_{i+1,j+1}(T_1)=D$, and then by (\ref{mr3_2}) we have
\begin{equation}\label{mr3_5}
 P_{i+1,j}(T_1) > P_{i+1,j+1}(T_1).
\end{equation}
At time $T_1$, for node $(i+1,j+1)$ we choose the control input as
 $\mathcal{C}_{i+1,j+1}(T_1)=(i+1,j)$, while for any other node we choose a neighbor which has the same strategy with it as the control input of $T_1$.
Then, under the $\mathcal{N}_C$-CEG, due to (\ref{mr3_3}) and (\ref{mr3_5}), we get that $S_{i+1,j+1}(T_1+1)=C$.
Meanwhile, the strategies of other nodes except $(i+1,j+1)$ keep unaltered at time $T_1+1$. We adopt the control input (\ref{mresult_1}) for the nodes in $\mathcal{N}\backslash \mathcal{N}_C$ at time $T_1+1$. Thus,
\begin{multline}\label{mr3_5_1}
P_{i,j+1}(T_1+1)\geq 3p_3+p_4>\max\big\{P_{i,j}(T_1+1),\\
P_{i-1,j+1}(T_1+1),P_{i+1,j+1}(T_1+1)\big\}.
\end{multline}
 By (\ref{mr3_5_1}) and (\ref{mresult_1_3}) we get
\begin{equation}\label{mr3_5_2}
S_{i,j+1}(T_1+2)=S_{i,j}(T_1+2)=S_{i+1,j+1}(T_1+2)=D.
\end{equation}
An evolution of nodes' strategies from (\ref{mr3_3}) to (\ref{mr3_5_2}) is shown in Fig. \ref{LFig5}.
Combining (\ref{mresult_1_2}), (\ref{mr3_3}), and (\ref{mr3_5_2}) yields
\begin{equation}\label{mr3_6}
n_{C}(T_1)> n_C(T_1+2).
\end{equation}
 \begin{figure*}
  \centering
  \includegraphics[width=6.5in]{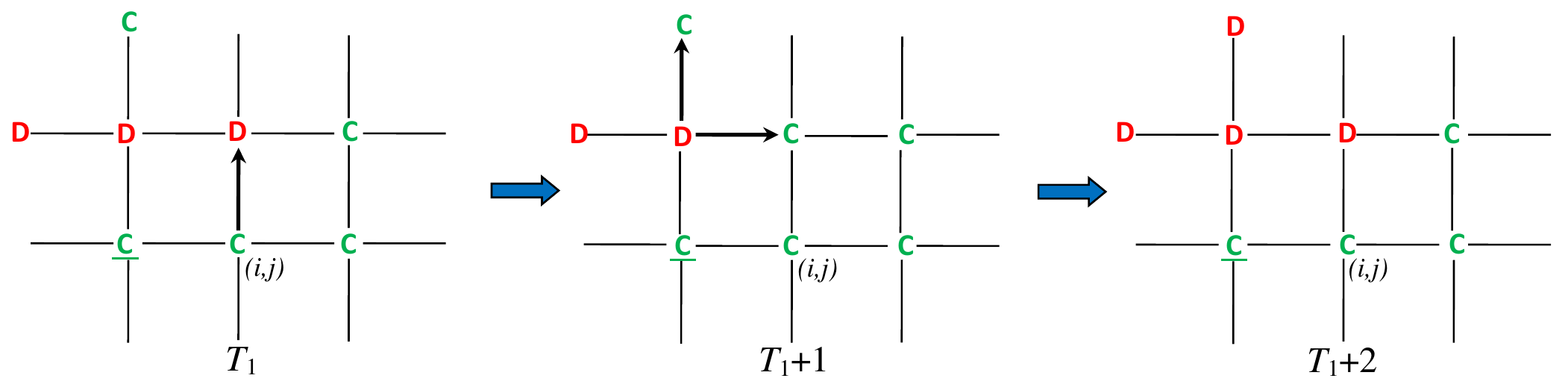}
  \caption{An evolution of nodes' strategies when $S_{i-1,j+1}(T_1)=D$ in Step 2, where ``$\underline{C}$" denotes the node belongs to  $\mathcal{N}_C$.}\label{LFig6}
\end{figure*}

If $S_{i-1,j+1}(T_1)=D$,  since $(i-1,j)\in \mathcal{N}_C$ and $\mathcal{N}_C$ is a $N_1\times M_1$ rectangular with $N_1\geq 2, M_1\geq 2$, we have $P_{i-1,j}(T_1)=3p_1+p_2$.
If $P_{i-1,j+1}(T_1)<P_{i-1,j}(T_1)$, with similar control inputs as Case II of the Step 3 in the proof of Theorem~\ref{Main_result} (see Fig. \ref{LFig1}), we can get (\ref{mr3_6}).
If $P_{i-1,j+1}(T_1)>P_{i-1,j}(T_1)$(Note that $P_{i-1,j+1}(T_1)\neq P_{i-1,j}(T_1)$ by the condition (\ref{cor1_0})), we have $P_{i-1,j+1}(T_1)\geq 2p_3+2p_4$, which indicates
\begin{equation}\label{mr3_7}
S_{i-2,j+1}(T_1)=C \mbox{~~or~~} S_{i-1,j+2}(T_1)=C.
\end{equation}
At time $T_1$, for node $(i,j+1)$ we choose the control input as
 $\mathcal{C}_{i,j+1}(T_1)=(i,j)$, while for any other node we choose a neighbor which has the same strategy with it as the control input of $T_1$. Then, under the $\mathcal{N}_C$-CEG,  it follows from (\ref{mr3_2})
 and (\ref{mr3_3}) that $S_{i,j+1}(T_1+1)=C$. Meanwhile, the strategies of other nodes except $(i,j+1)$ keep unvaried at time $T_1+1$.
 Thus, we have
\begin{equation}\label{mr3_8}
S_{i-1,j+1}(T_1+1)=D ~\mbox{and}~ P_{i-1,j+1}(T_1+1)\geq 3p_3+p_4.
\end{equation}
Since $(i-1,j)\in\mathcal{N}_C$, $(i-1,j+1)\notin\mathcal{N}_C$, and $M_1\notin\{M-1,M-2\}$, we have
$(i-2,j+1)\notin\mathcal{N}_C$ and $(i-1,j+2)\notin\mathcal{N}_C$.
We adopt the control input (\ref{mresult_1}) at time $T_1+1$. By (\ref{mr3_8}) and (\ref{mresult_1_3}) we get
\begin{equation}\label{mr3_9}
S_{i,j+1}(T_1+2)=S_{i-2,j+1}(T_1+2)=S_{i-1,j+2}(T_1+2)=D.
\end{equation}
An evolution of nodes' strategies for the above process is shown in Fig. \ref{LFig6}.
Combining (\ref{mresult_1_2}), (\ref{mr3_7}), and (\ref{mr3_9}) yields (\ref{mr3_6}).

Furthermore, when $(i+1,j)\in \mathcal{N}_C$, or $(i,j-1)\in\mathcal{N}_C$, or $\{(i-1,j), (i,j-1), (i+1,j)\}\cap \mathcal{N}_C =\emptyset$, with the similar discussion to the case $(i-1,j)\in \mathcal{N}_C$
we get (\ref{mr3_6}).

For $t\geq T_1+2$  we repeatedly carry out Step 1 and the above process until  there is no cooperation node in $\mathcal{N}\backslash\mathcal{N}_C$ which has $3$ cooperation neighbors, and the strategies of all nodes keep unchanged under control input (\ref{mresult_1}). We record the stop time as $T_2$.
If $S(T_2)=S^*$, our assertion holds. Otherwise, we perform the following Step 3.

\textbf{Step 3}:
If there is a cooperation node $(i,j)\notin \mathcal{N}_C$  which has two cooperation neighbors at time $T_2$,  because $\mathcal{N}_C$ forms a $N_1\times M_1$ rectangular area with $N_1\in [2,N-3]\cup\{N\}, M_1\in [2,M-3]\cup\{M\}$,  the node $(i,j)$ has at least one cooperation neighbor which does not belong to $\mathcal{N}_C$. Without loss of generality we assume
\begin{equation}\label{mr3_10}
S_{i,j}(T_2)=S_{i+1,j}(T_2)=C, (i+1,j)\notin\mathcal{N}_C,  S_{i,j+1}(T_2)=D.
\end{equation}
It is derived from (\ref{mr3_10}) and (\ref{mresult_1_3}) that
\begin{equation}\label{mr3_10_1}
P_{i,j+1}(T_2)\leq P_{i,j}(T_2)=2p_1+2p_2.
\end{equation}
Also, by (\ref{mr3_10}) and the Step 2 we have
\begin{equation}\label{mr3_10_2}
 P_{i+1,j}(T_2)\leq 2p_1+2p_2.
\end{equation}
One the other hand, by the conditions of $p_1,\ldots,p_4$ we have $p_1+p_2<p_3+p_4$, so by (\ref{mr3_10_1})
at time $T_2$ the defection node $(i,j+1)$ has only one cooperation neighbor $(i,j)$, which indicates
 $S_{i+1,j+1}(T_2)=D$. Then, by (\ref{mresult_1_3}) and  the condition (\ref{cor1_0}) we have
\begin{equation}\label{mr3_10_3}
P_{i+1,j+1}(T_2)< P_{i+1,j}(T_2).
\end{equation}
Combining (\ref{mr3_10_2}), (\ref{mr3_10_3}) and the fact $p_1+p_2<p_3+p_4$ yields that
 the node $(i+1,j+1)$ has only one cooperation neighbor at time $T_2$.
We choose $\mathcal{C}_{i+1,j+1}(T_2)=(i+1,j)$ for node $(i+1,j+1)$, while for any other node we choose a neighbor which has the same strategy with it as the control input.  Then, by the $\mathcal{N}_C$-CEG, at time $T_2+1$ the strategy of the node $(i+1,j+1)$ changes from $D$ to $C$, while the other nodes keep strategy unaltered. Thus,
we have
\begin{multline}\label{mr3_10_4}
P_{i,j+1}(T_2+1)=2p_3+2p_4\\
>2p_1+2p_2=P_{i,j}(T_2+1)>P_{i+1,j+1}(T_2+1).
\end{multline}
We adopt the control input (\ref{mresult_1}) at time $T_2+1$. According to (\ref{mr3_10_4})
and (\ref{mresult_1_3}) we can get
\begin{equation}\label{mr3_10_5}
S_{i,j+1}(T_2+2)=S_{i,j}(T_2+2)=S_{i+1,j+1}(T_2+2)=D.
\end{equation}
The evolution of nodes' strategies from (\ref{mr3_10}) to (\ref{mr3_10_5}) is shown in Fig. \ref{LFig2}
but using $T_2$ instead of $T_1$.
Combining (\ref{mr3_10_5}) with (\ref{mresult_1_2}) and (\ref{mr3_10}) yields
\begin{equation}\label{mr3_11}
n_{C}(T_2)> n_C(T_2+2).
\end{equation}

For $t\geq T_2+2$  we repeatedly carry out Steps 1-2 and the above process until  there is no cooperation node in $\mathcal{N}\backslash\mathcal{N}_C$ with $2$  or $3$ cooperation neighbors, and the strategies of all nodes keep unchanged under control input (\ref{mresult_1}). We record the stop time as $T_3$.
If $S(T_3)=S^*$, our assertion holds. Otherwise, we perform the following Step 4.

\textbf{Step 4}: If there exists a cooperation node $(i,j)\notin \mathcal{N}_C$  which has exact one cooperation neighbor at time $T_3$, without loss of generality we assume $S_{i+1,j}(T_3)=C$.

If $(i+1,j)\in\mathcal{N}_C$, because $\mathcal{N}_C$ forms a rectangular area, without loss of generality we can assume $(i+1,j+1)\in\mathcal{N}_C$. Thus, by the fact of $S_{i,j+1}(T_3)=D$ and the conditions of $p_1,\ldots,p_4$ we have
\begin{eqnarray*}\label{mr3_12}
P_{i,j+1}(T_3)\geq 2p_3+2p_4>p_1+3p_2=P_{i,j}(T_3),
\end{eqnarray*}
which is conflict with (\ref{mresult_1_3}). According to the definition of $T_3$ the case that $(i+1,j)\in\mathcal{N}_C$ will not happen.

If $(i+1,j)\notin\mathcal{N}_C$, by the similar method to the Case II of the Step 2 in the proof of Theorem
\ref{Main_result2} we can obtain
\begin{equation}\label{mr3_13}
n_{C}(T_3)> n_C(T_3+2).
\end{equation}

For $t\geq T_3+2$  we repeatedly carry out Steps 1-3 and the above process until $\{S^*\}$ is reached. Let $T_4$ be the stop time. By (\ref{mresult_2}), (\ref{mr3_6}), (\ref{mr3_11}), and (\ref{mr3_13}) we have
\begin{equation}\label{mr3_12}
 T_4\leq 2(NM-N_1M_1-1).
\end{equation}

Next we show  that under the CEG,  $\{S_C\}$ is finite-time reachable from
$\{S^*\}$. We continue to carry out the following Step 5:

\textbf{Step 5}: Without loss of generality we assume $N_1\in [2,N-3]$.
We first consider the case when $N-N_1$ is an even number.  At time $T_4$, for nodes $(i_0+N_1,j), j_0\leq j\leq j_0+M_1-1$, we choose their control inputs $\mathcal{C}_{i_0+N_1,j}(T_4)=(i_0+N_1-1,j)$; for
 nodes $(i_0-1,j), j_0\leq j\leq j_0+M_1-1$, we choose their control inputs $\mathcal{C}_{i_0-1,j}(T_4)=(i_0,j)$; for other nodes in $\mathcal{N}\backslash\mathcal{N}_C$   we choose defection neighbors as their control inputs.
 Because $S(T_4)=S^*$,  by (\ref{mr3_1})  and the fact that $\mathcal{N}_C$ is a rectangular with $N-N_1>1,$ and $M-M_1=0$ or $M-M_1>1$, if $(i,j)\notin\mathcal{N}_C$ and $(i,j)$ has a cooperation neighbor then $P_{i,j}(T_4)=p_3+3 p_4<2p_1+2p_2$. Thus, under the $\mathcal{N}_C$-CEG, we have  $S_{i_0+N_1,j}(T_4+1)=S_{i_0-1,j}(T_4+1)=C$  for $j_0\leq j\leq j_0+M_1-1$, while the strategies of other nodes keep unvaried. The cooperation nodes form a $(N_1+2)\times M_1$ rectangular at time $T_4+1$. With the similar process as above we can get that the cooperation nodes form a $N\times M_1$ rectangular (torus)  at time $T_4+(N-N_1)/2$.

If $N-N_1$ is an odd number,  at time $T_4$, for nodes $(i_0+N_1,j), j_0\leq j\leq j_0+M_1-1$, we choose $\mathcal{C}_{i_0+N_1,j}(T_4)=(i_0+N_1-1,j)$, while for other nodes in $\mathcal{N}\backslash\mathcal{N}_C$   we choose defection neighbors as their control inputs. Then, the cooperation nodes form a $(N_1+1)\times M_1$ rectangular at time $T_4+1$. Because $N-N_1-1>1$ is an even number, by the similar process as the case when $N-N_1$ is even, we can get the cooperation nodes form a $N\times M_1$ rectangular (torus)  at time $T_4+(N-N_1+1)/2$.

Let $T_5:=T_4+\lceil(N-N_1)/2\rceil$. With the similar method as above we can choose appropriate control inputs such that the cooperation nodes form a $N\times M$ rectangular (torus)  at time $T_5+\lceil(M-M_1)/2\rceil$, which means $S(T_5+\lceil(M-M_1)/2\rceil)=S_C$.

\renewcommand{\thesection}{\Roman{section}}
\section{Simulations}\label{Simulations}
\renewcommand{\thesection}{\arabic{section}}

This section tries to find the sufficient and necessary condition of convergence of SEG from simulations. All simulations are carried out by assuming that $M=N=10$,
the initial strategy of each node is chosen to be $C$ or $D$ with equal probability,  and the imitation probability of each node equals $1$ if its payoff is less than the payoff of its neighbor,
i.e.,
\begin{eqnarray*}\label{updaterule_1}
&&\Prob \left\{S_{i,j}(t+1)=S(\mathcal{R}_{i,j}(t))| P_{i,j}(t)<P(\mathcal{R}_{i,j}(t)) \right\}\\
&&~=1, ~~~~\forall 1\leq i\leq N, 1\leq j\leq M, t\geq 0,
\end{eqnarray*}
where $\mathcal{R}_{i,j}(t)$ is independently and uniformly selected from $(i,j)$'s  neighbor nodes.

\begin{figure}
  \centering
  \subfigure[$c=0.74$, $t=0$]{
    \includegraphics[scale=0.27]{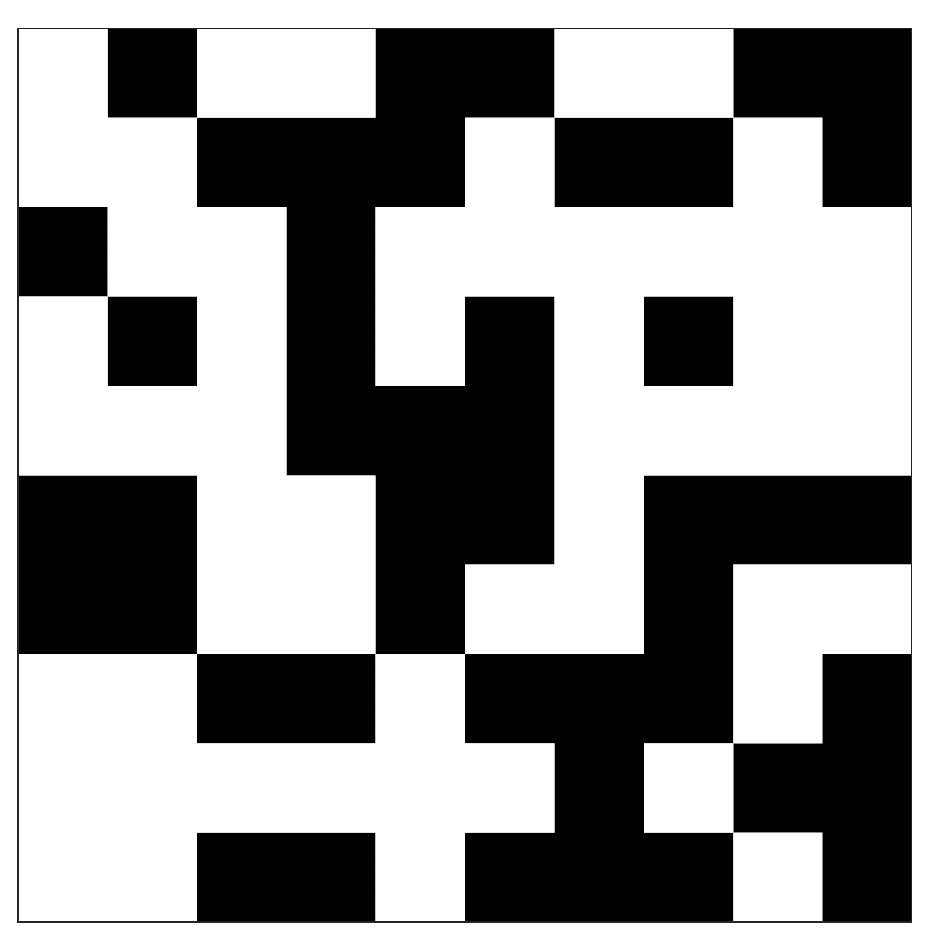}}
  \hspace{0.02in} 
  \subfigure[$c=0.74$, $t=100$]{
    \includegraphics[scale=0.27]{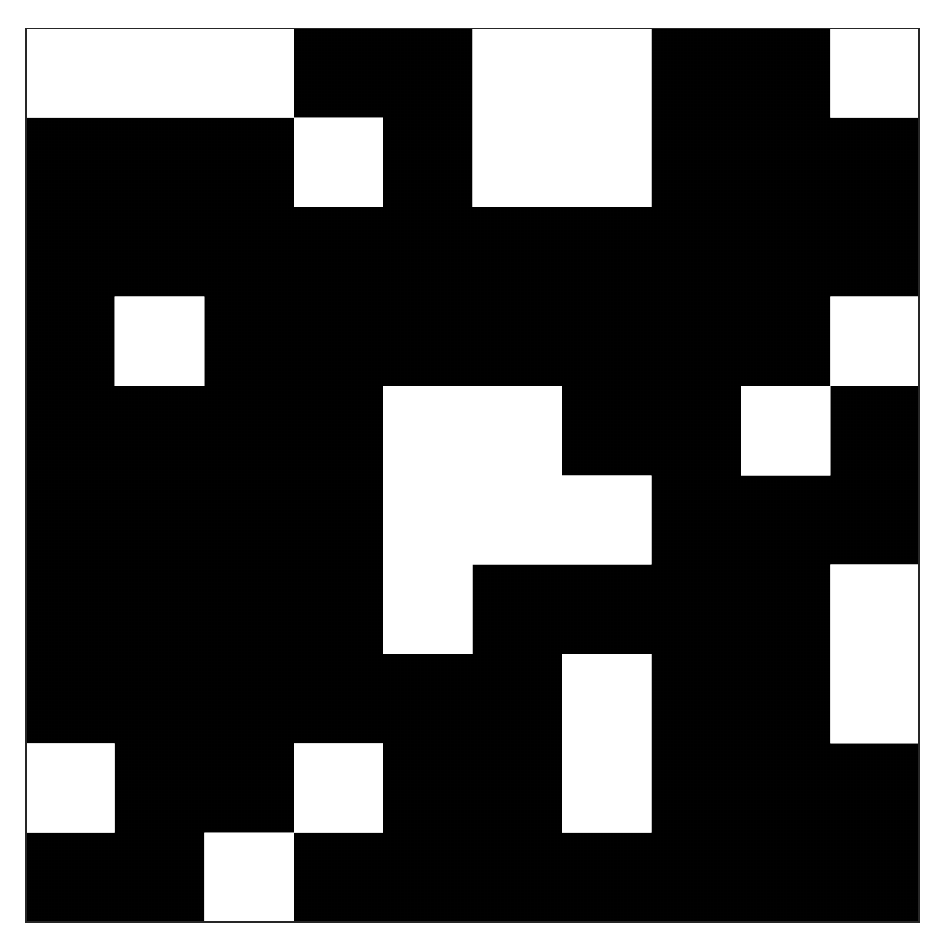}}
    \hspace{0.02in} 
  \subfigure[$c=0.74$, $t=10^7$]{
    \includegraphics[scale=0.27]{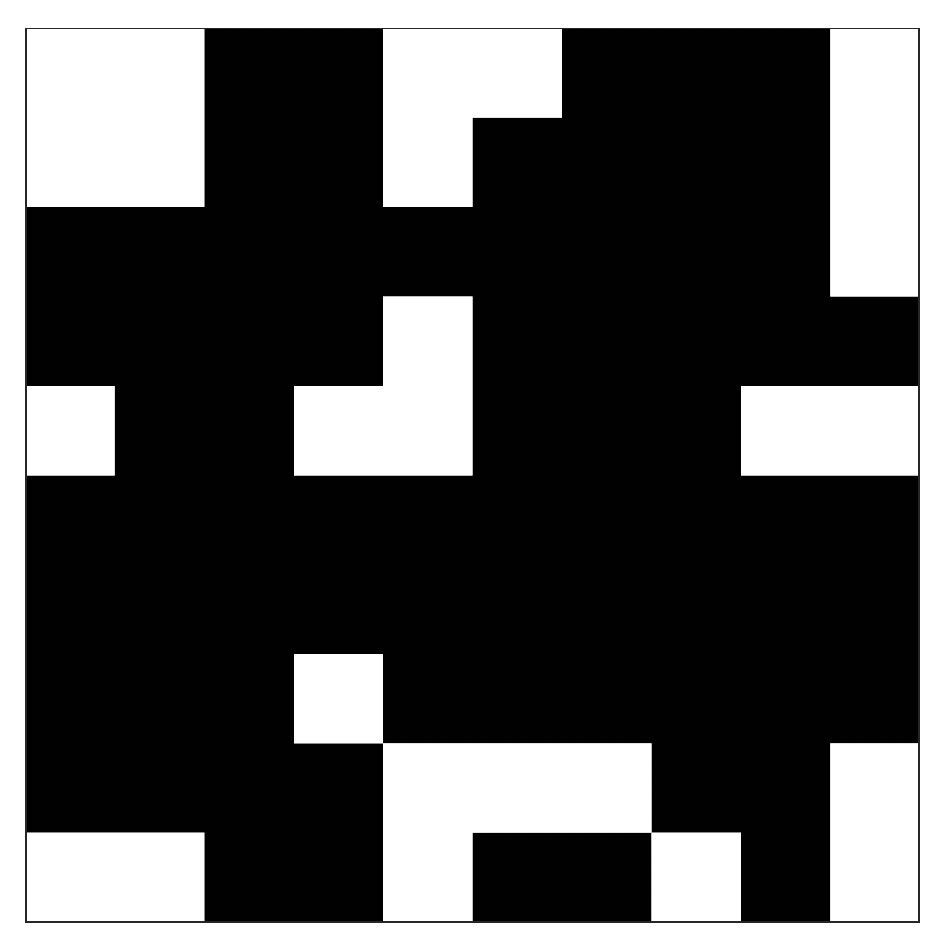}}
  \subfigure[$c=0.75, t=0$]{
    \includegraphics[scale=0.27]{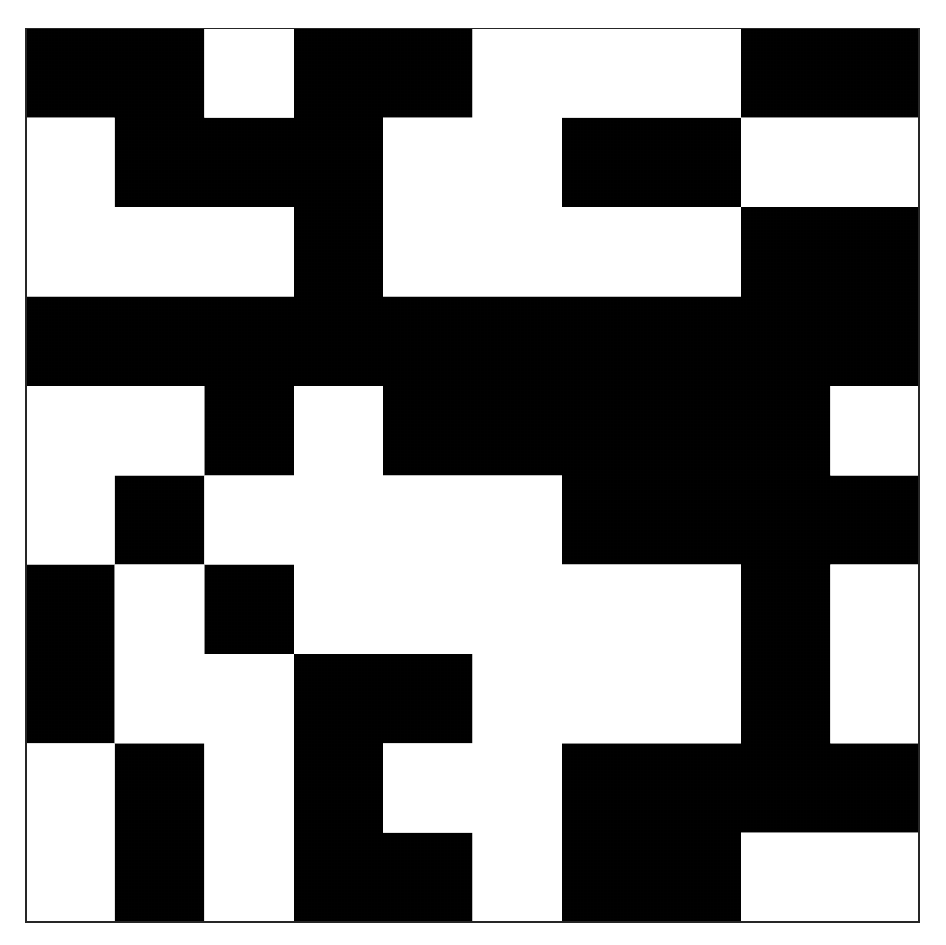}}
  \hspace{0.02in} 
  \subfigure[$c=0.75, t=100$]{
    \includegraphics[scale=0.27]{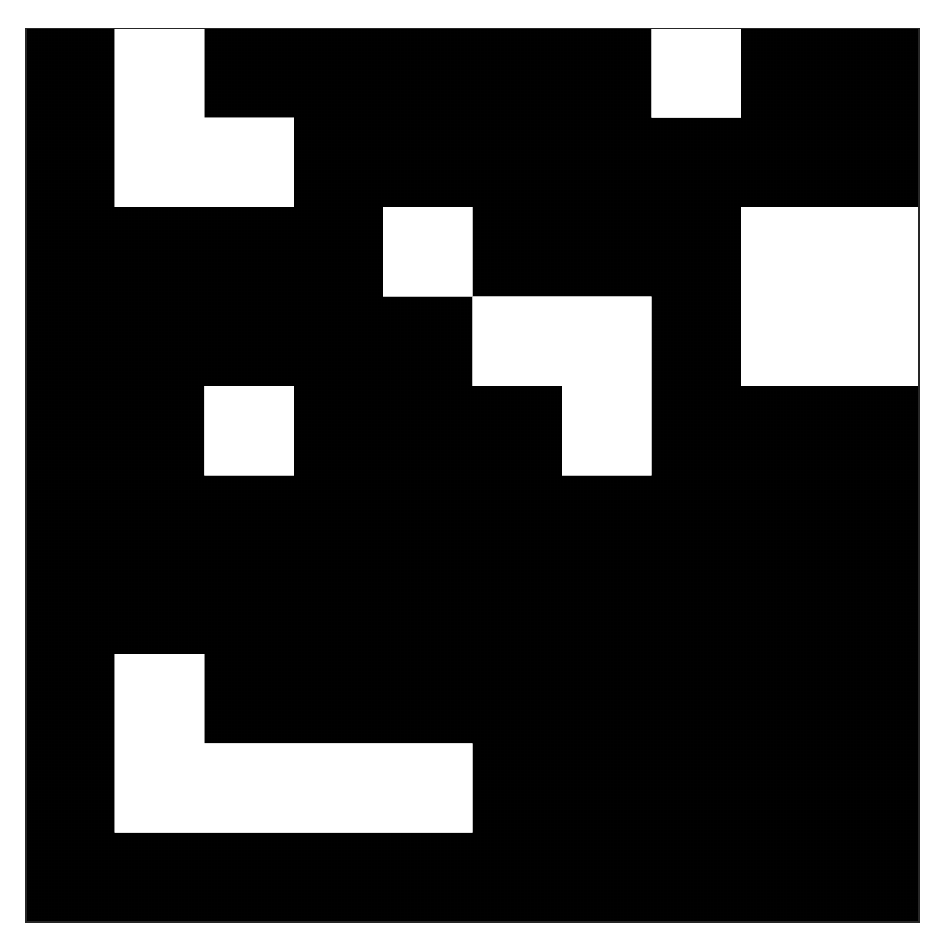}}
    \hspace{0.02in} 
  \subfigure[$c=0.75, t=111$]{
    \includegraphics[scale=0.27]{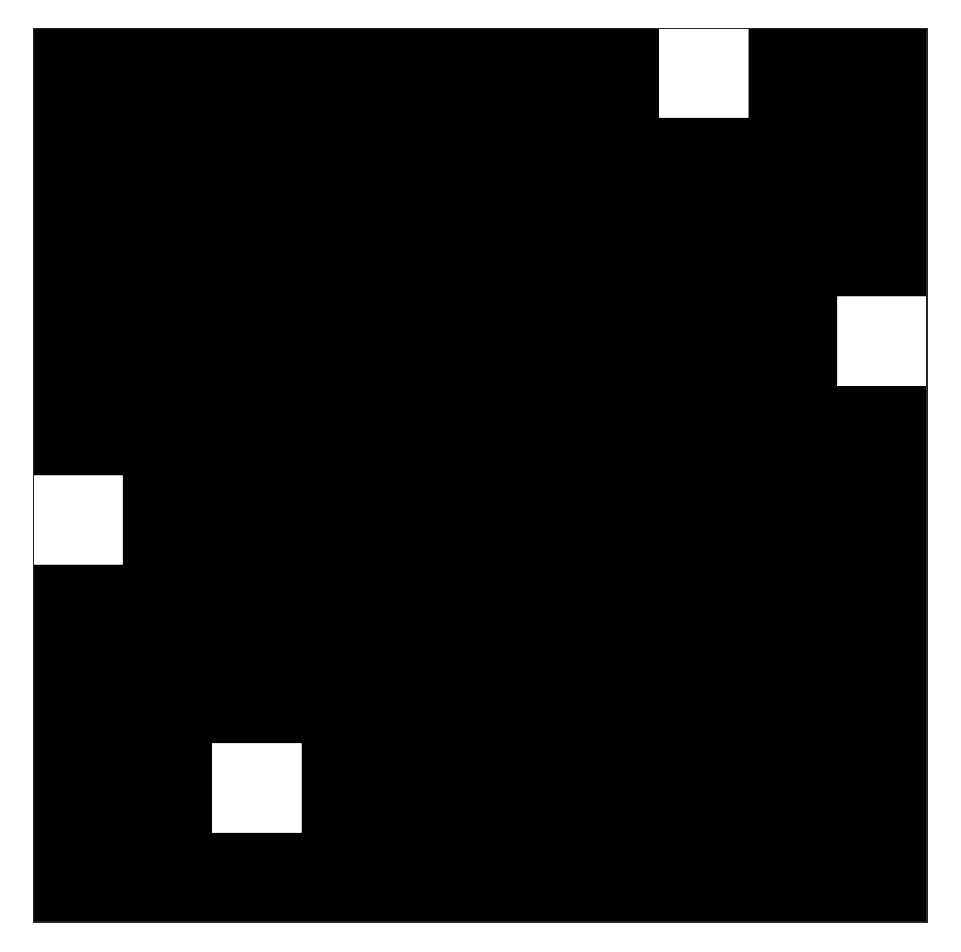}}
  \subfigure[$c=0.76, t=0$]{
    \includegraphics[scale=0.27]{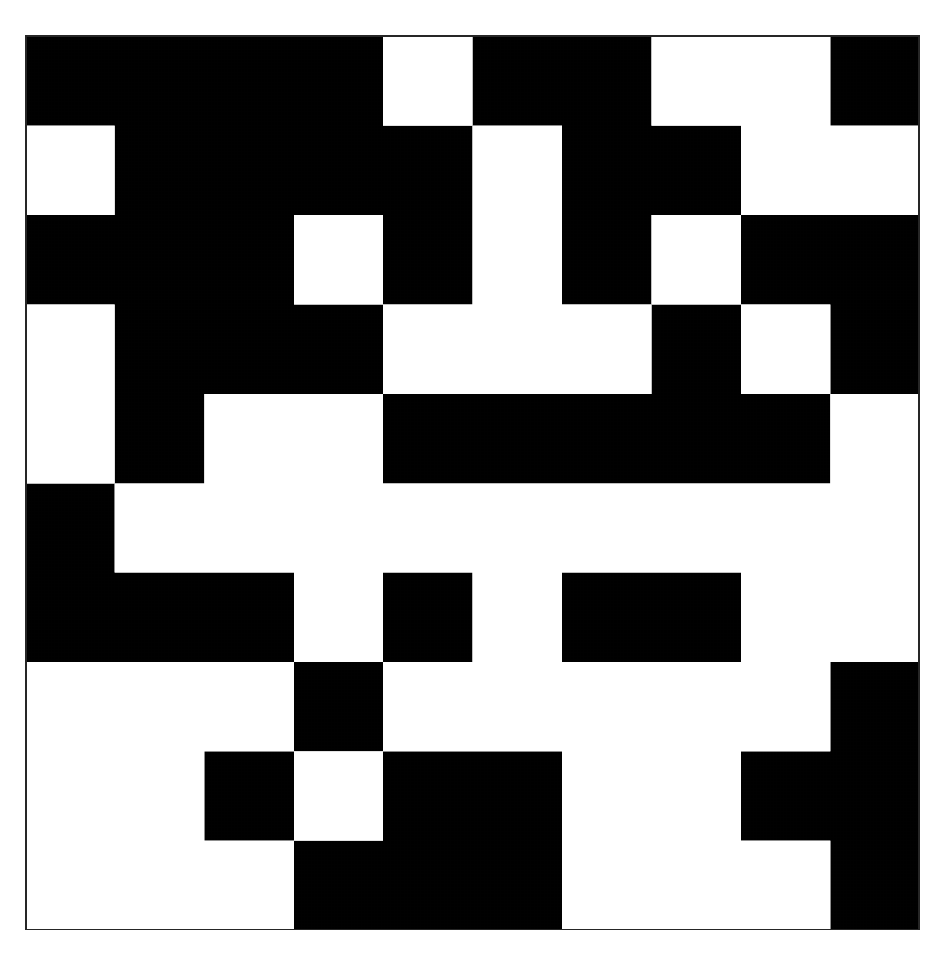}}
  \hspace{0.02in} 
  \subfigure[$c=0.76, t=100$]{
    \includegraphics[scale=0.27]{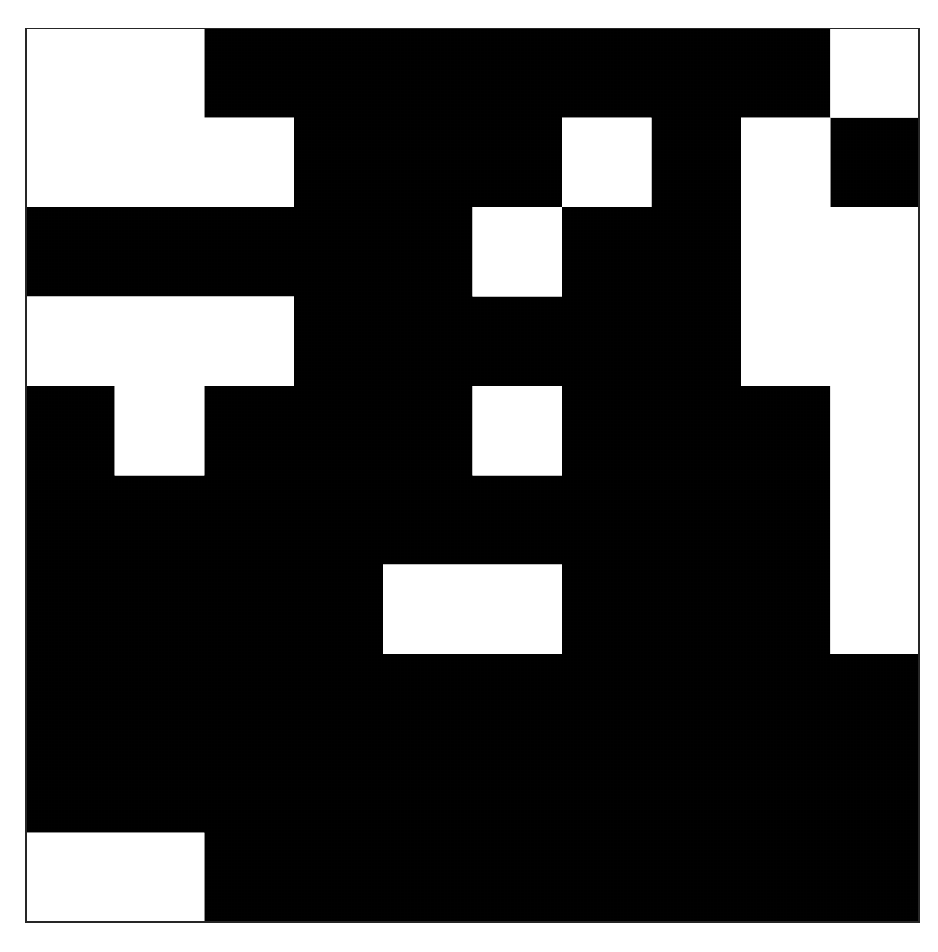}}
    \hspace{0.02in} 
  \subfigure[$c=0.76, t=288$]{
    \includegraphics[scale=0.27]{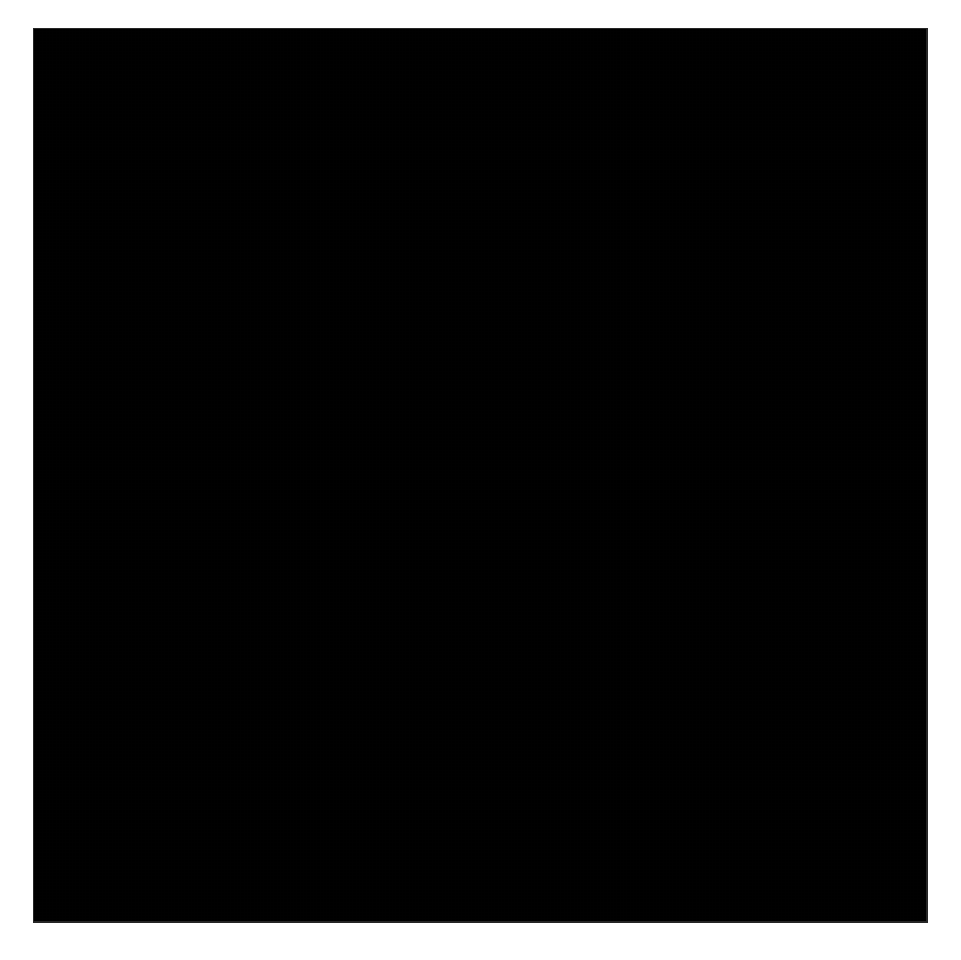}}
  \caption{The evolutionary snowdrift game with $c=0.74, 0.75, 0.76$ respectively. The system converges to a fixed state at time $t=111$ for $c=0.75$, and at time $t=288$ for $c=0.76$, however  does not converge   at time $t=10^7$ for $c=0.74$.}
  \label{esg76}
\end{figure}

First, we consider the evolutionary snowdrift game with the classic payoff matrix
\begin{table}[hbp]
\centering
\begin{tabular}{| c | c | c |}
\hline
      & C & D \\
\hline
      C & $1-c/2$ & $1-c$  \\
\hline
      D & $1$ & $0$ \\
 \hline
\end{tabular}
\end{table}\\
satisfying $0<c<1$. By Theorem \ref{Main_result2} if $c>3/4$ then the evolutionary snowdrift game converges to a fixed state a.s. in finite time.
Simulations show that the evolutionary snowdrift game does not converge to a fixed state when $c=0.74$, however converges fast
  when  $c=0.75,0.76$ (see Fig. \ref{esg76}. Throughout this section the white square denotes ``C", while the black square denotes ``D"), which imply that $3/4$ may be the critical value of $c$ for the convergence of evolutionary snowdrift game.

\begin{figure}
  \centering
  \subfigure[$b=0.59$, $t=0$]{
    \includegraphics[scale=0.27]{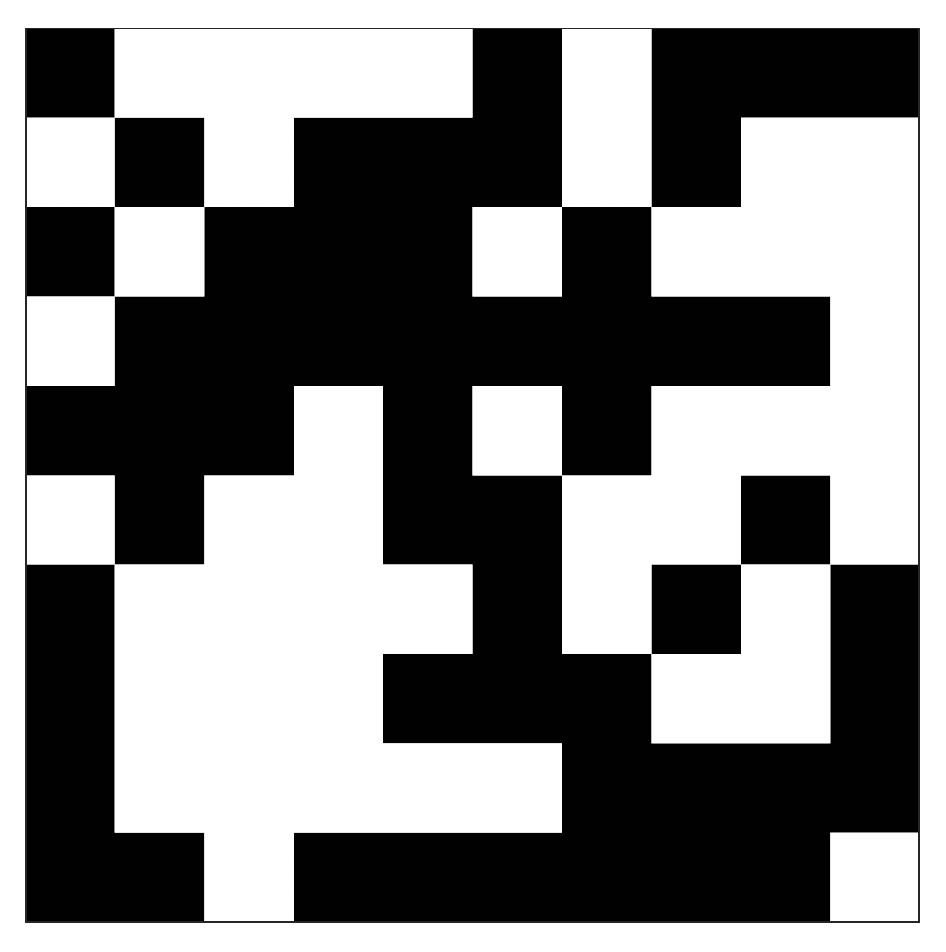}}
  \hspace{0.02in} 
  \subfigure[$b=0.59$, $t=100$]{
    \includegraphics[scale=0.27]{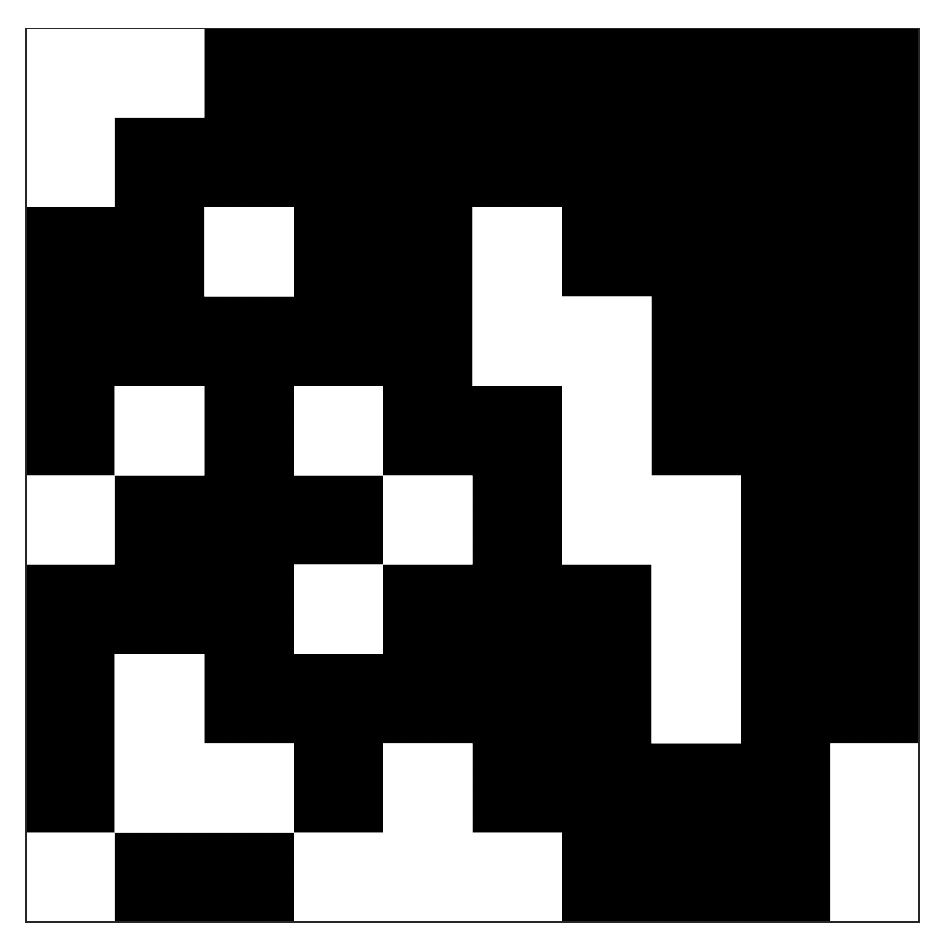}}
    \hspace{0.02in} 
  \subfigure[$b=0.59$, $t=10^7$]{
    \includegraphics[scale=0.27]{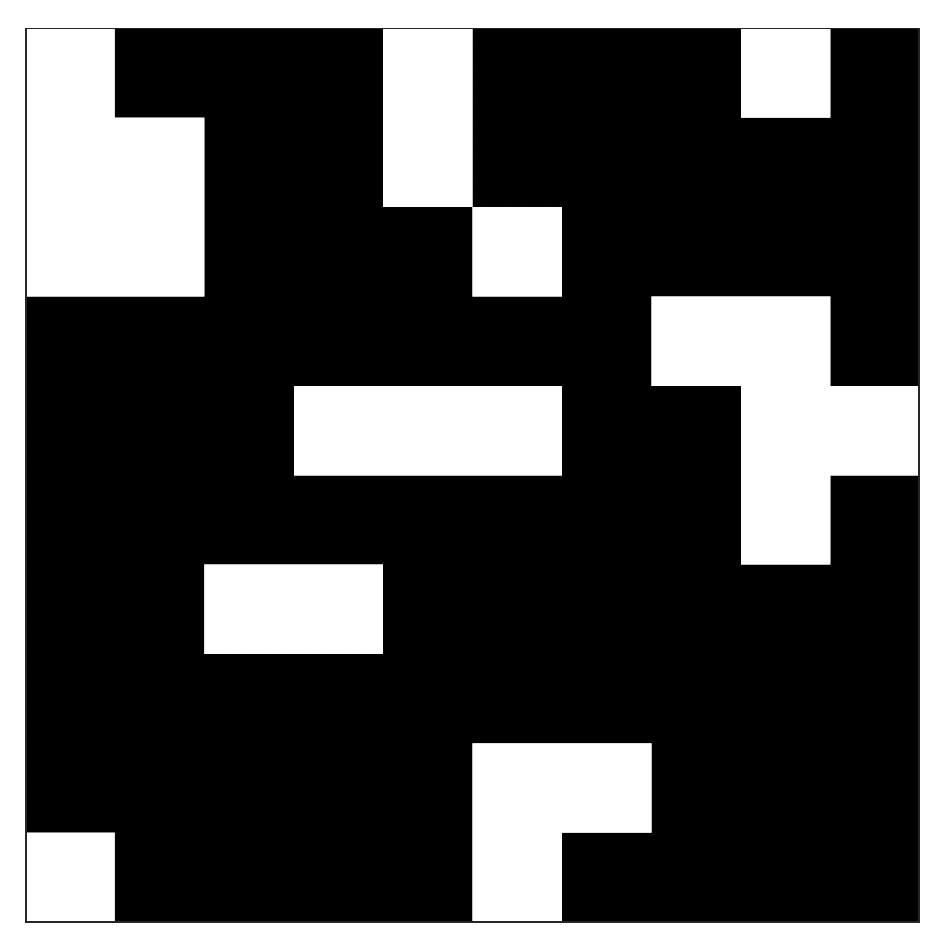}}
  \subfigure[$b=0.6, t=0$]{
    \includegraphics[scale=0.27]{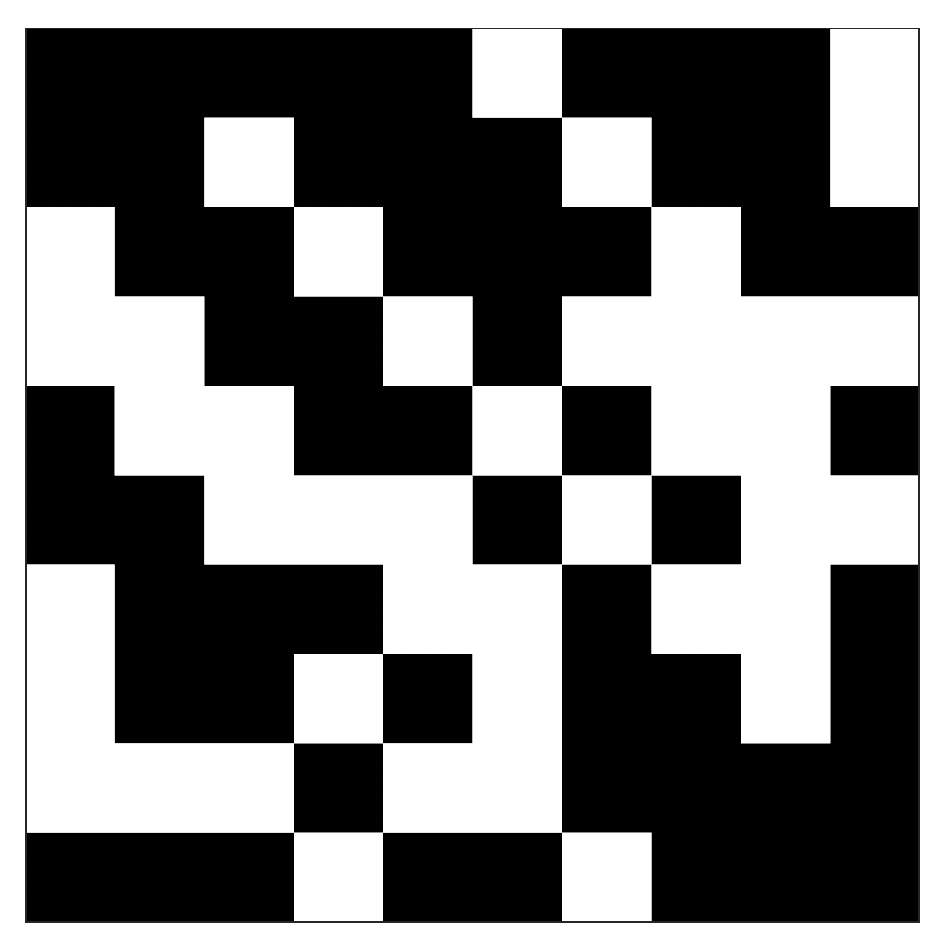}}
  \hspace{0.02in} 
  \subfigure[$b=0.6, t=100$]{
    \includegraphics[scale=0.27]{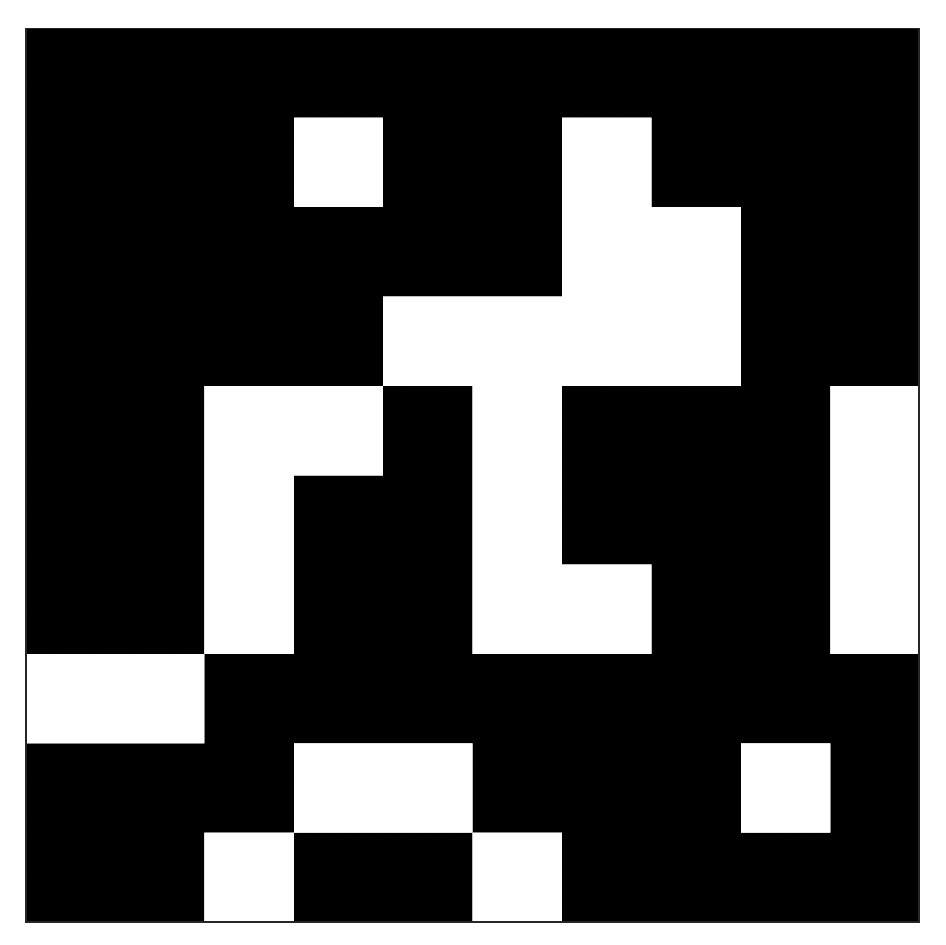}}
    \hspace{0.02in} 
  \subfigure[$b=0.6, t=10^7$]{
    \includegraphics[scale=0.27]{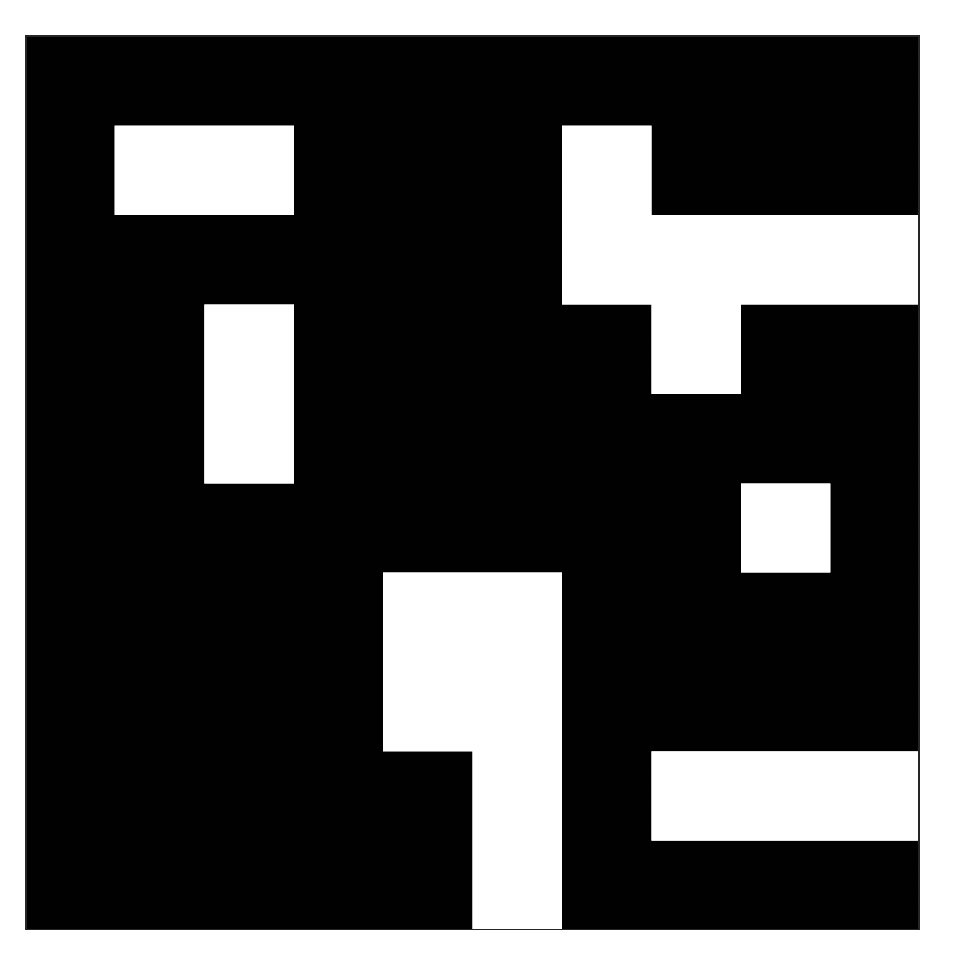}}
  \subfigure[$b=0.61, t=0$]{
    \includegraphics[scale=0.27]{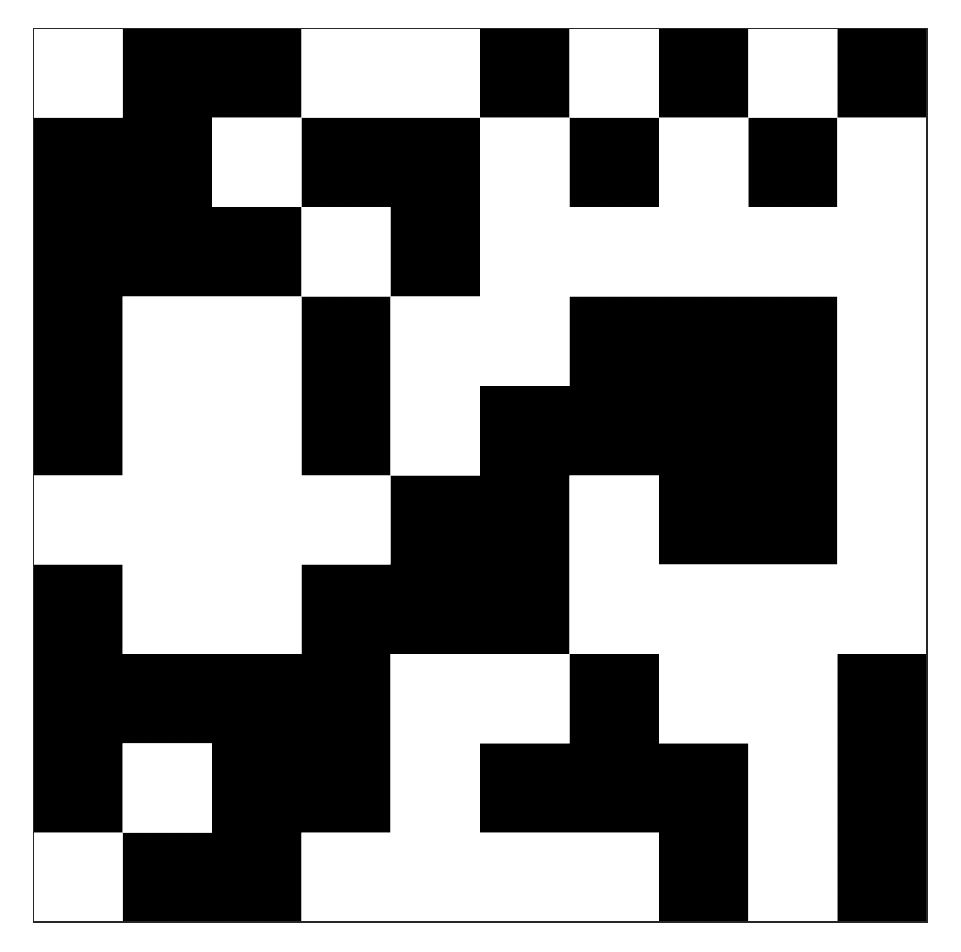}}
  \hspace{0.02in} 
  \subfigure[$b=0.61, t=100$]{
    \includegraphics[scale=0.27]{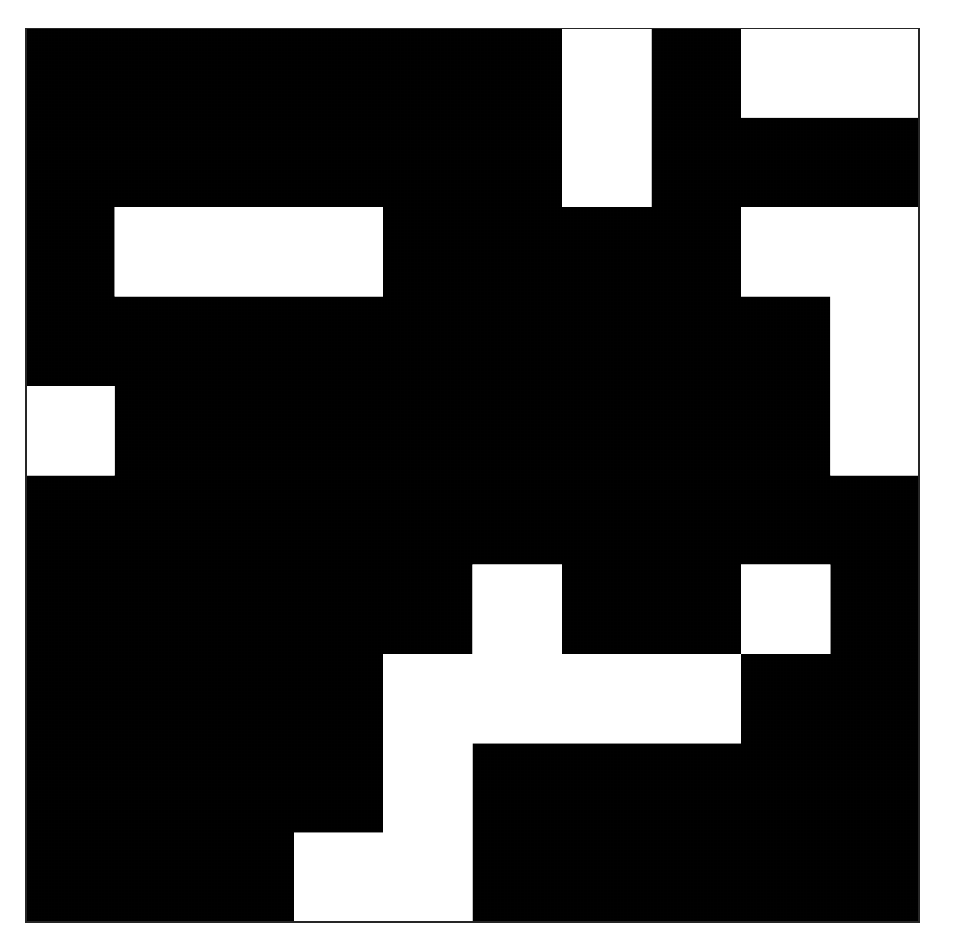}}
    \hspace{0.02in} 
  \subfigure[$b=0.61, t=176$]{
    \includegraphics[scale=0.27]{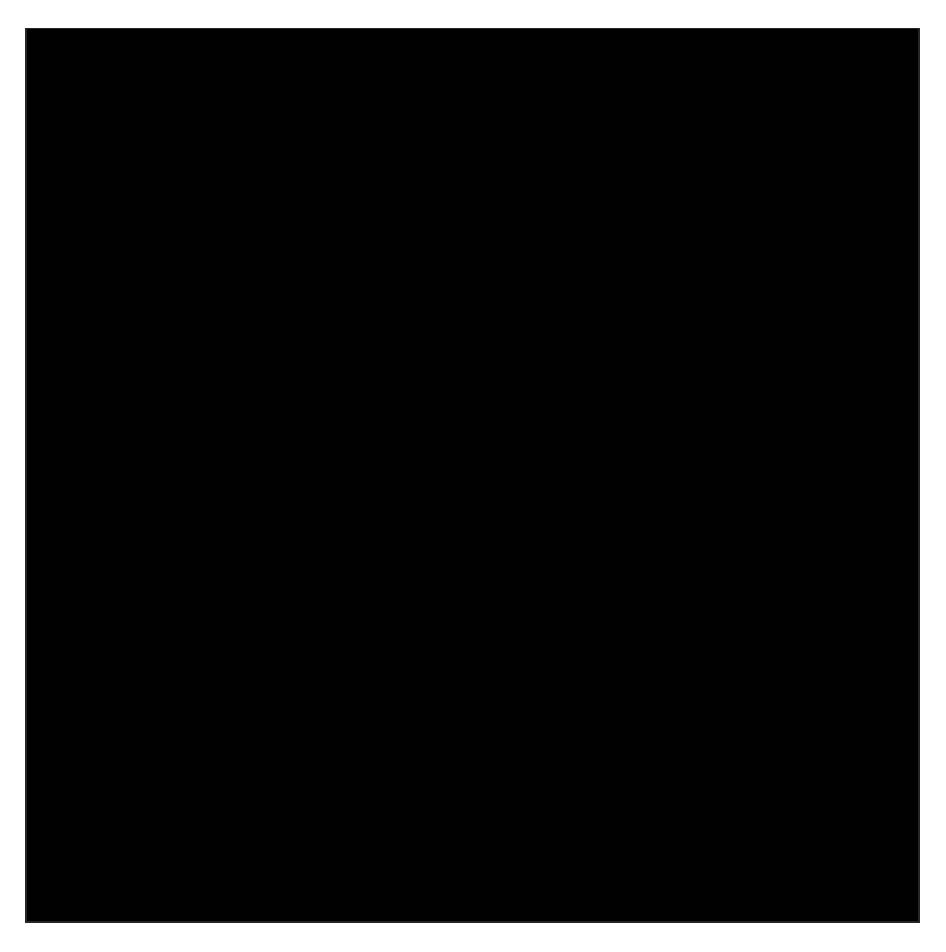}}
  \caption{The evolutionary hawk-dove game with $b=0.59, 0.6, 0.61$ respectively. The system converges to a fixed state at time $t=176$ for $b=0.61$, however  does not converge  at time $t=10^7$ for $b=0.59, 0.6$.}
  \label{ehd}
\end{figure}

Second, we simulate the evolutionary  hawk-dove game with the classic payoff matrix
\begin{table}[!hbp]
\centering
\begin{tabular}{| c | c | c |}
\hline
      & C & D \\
\hline
      C & $b/2$ & $0$  \\
\hline
      D & $b$ & $(b-1)/2$ \\
 \hline
\end{tabular}
\end{table}\\
satisfying $0<b<1$. By Theorem \ref{Main_result2} if $b>3/5$ then the evolutionary hawk-dove game converges to a fixed state a.s. in finite time.
Simulations show that the evolutionary hawk-dove game does not converge to a fixed state when $b=0.59, 0.6$, however converges fast when
   $b=0.61$ (see Fig. \ref{ehd}), which imply that $3/5$ may be the critical value of $b$ for the convergence of evolutionary hawk-dove game.

Third, we simulate the  evolutionary  chicken game with the classic payoff matrix
\begin{table}[hbp]
\centering
\begin{tabular}{| c | c | c |}
\hline
      & C & D \\
\hline
      C & $b/2$ & $0$  \\
\hline
      D & $b$ & $-1$ \\
 \hline
\end{tabular}
\end{table}\\
satisfying $b>0$. By Theorem \ref{Main_result2} if $b>3$ then the evolutionary chicken game converges to a fixed state a.s. in finite time.
Simulations show that the evolutionary chicken game does not converge to a fixed state when $b=2.9$, however converges fast
  when  $b=3,3.1$ (The figures are similar to Fig. \ref{esg76}, then to save space we do not put them here), which imply that $3$ may be the critical value of $b$ for the convergence of evolutionary  chicken game.

\begin{figure}
  \centering
  \subfigure[$r=0.32$, $t=0$]{
    \includegraphics[scale=0.27]{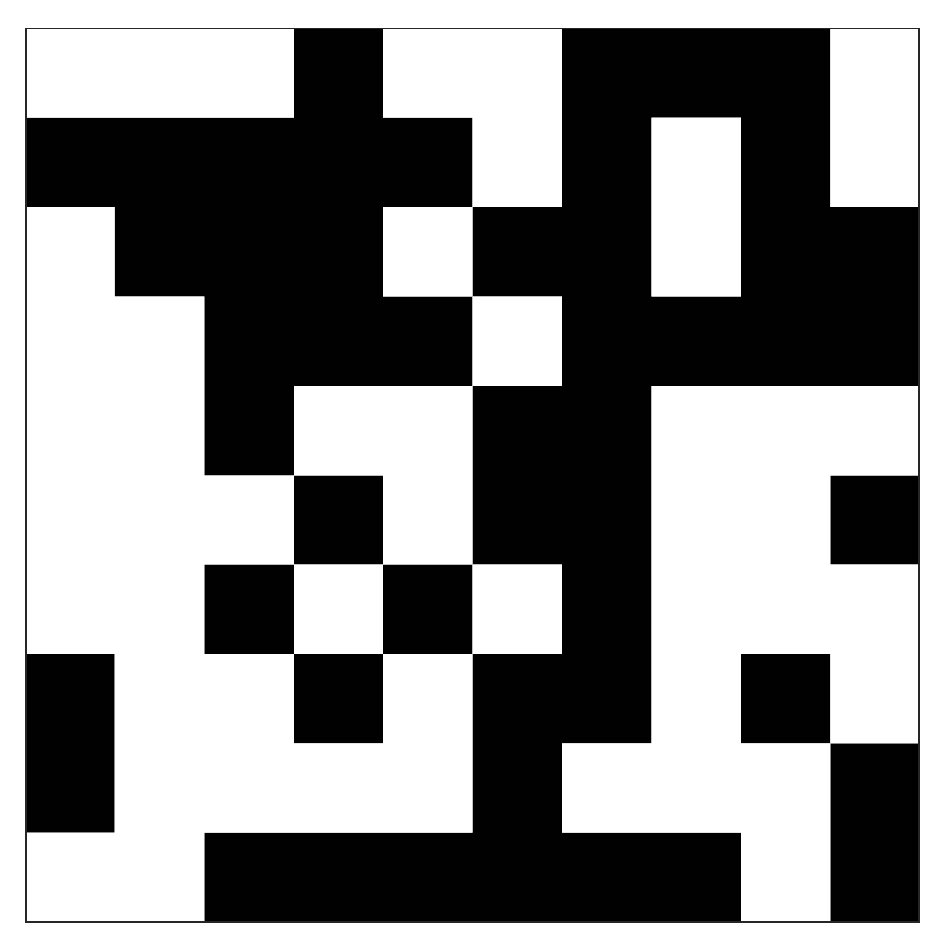}}
  \hspace{0.02in} 
  \subfigure[$r=0.32$, $t=5$]{
    \includegraphics[scale=0.27]{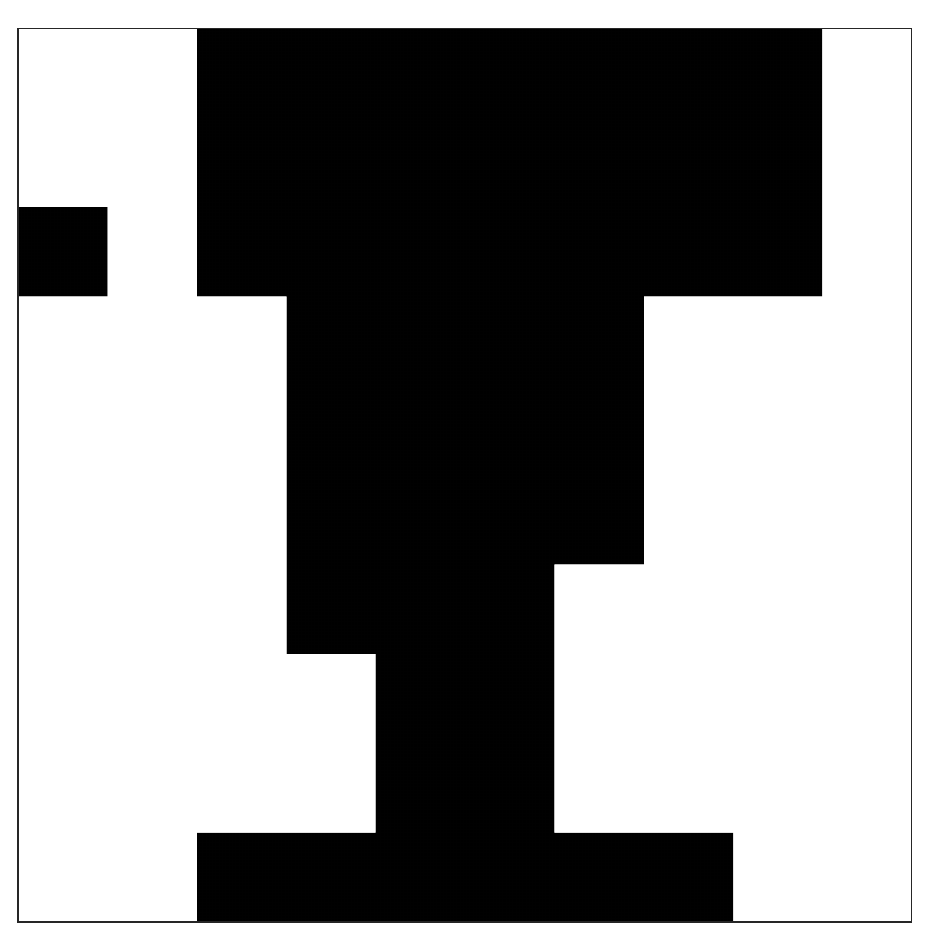}}
    \hspace{0.02in} 
  \subfigure[$r=0.32$, $t=13$]{
    \includegraphics[scale=0.27]{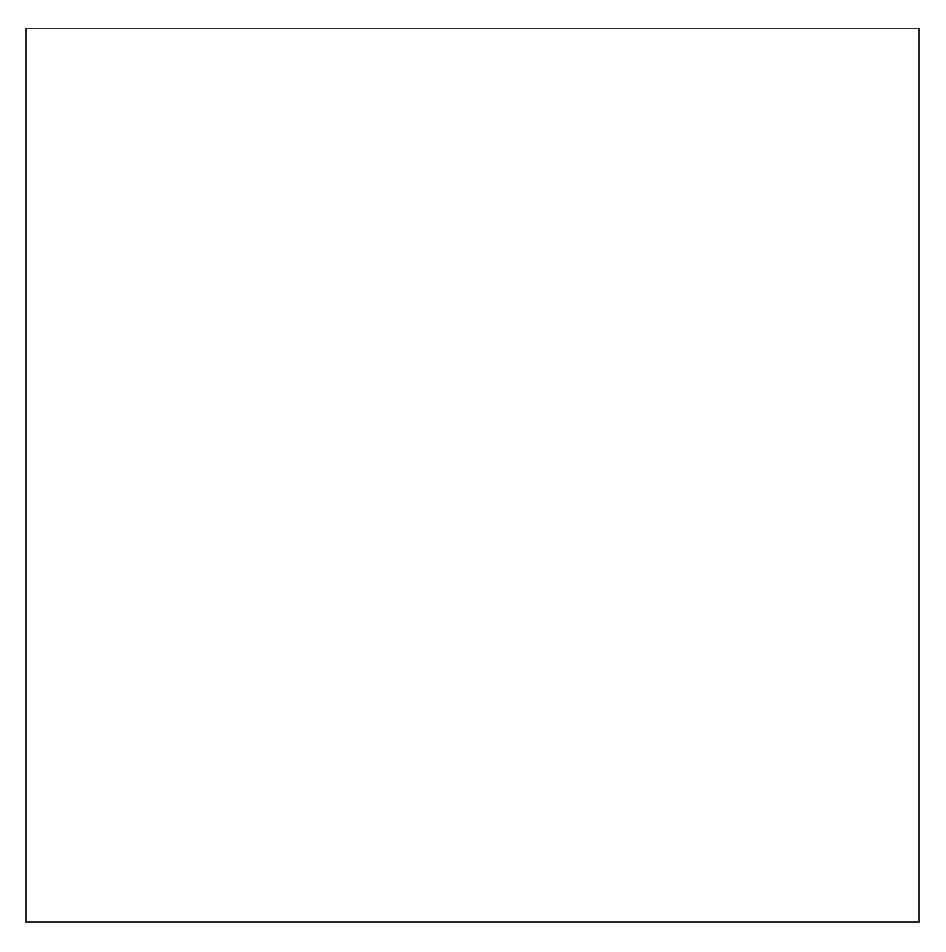}}
  \subfigure[$r=1/3, t=0$]{
    \includegraphics[scale=0.27]{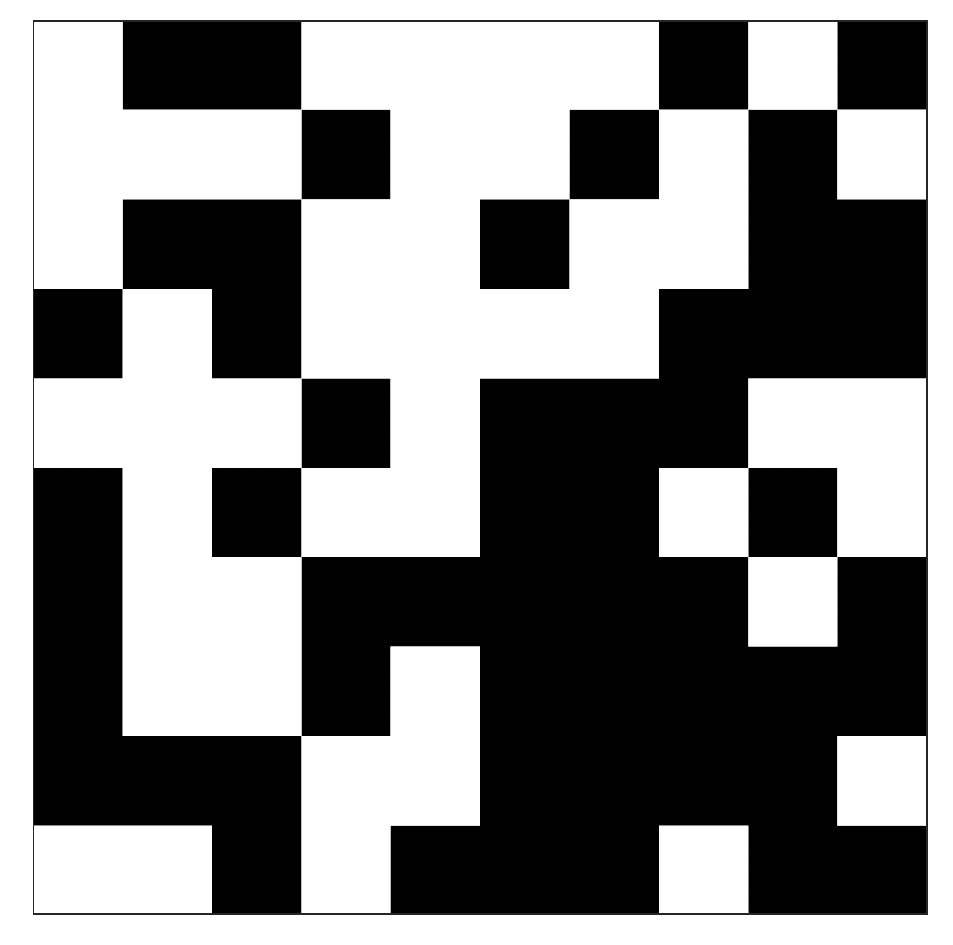}}
  \hspace{0.02in} 
  \subfigure[$r=1/3, t=5$]{
    \includegraphics[scale=0.27]{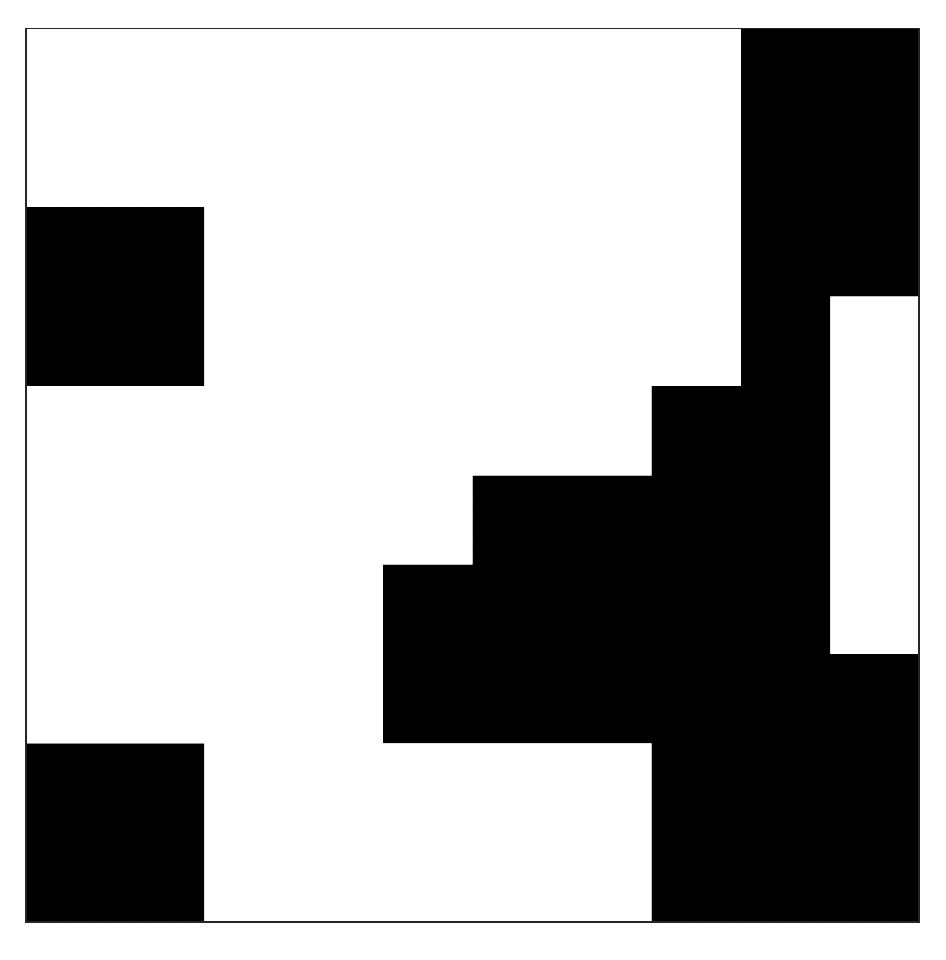}}
    \hspace{0.02in} 
  \subfigure[$r=1/3, t=11$]{
    \includegraphics[scale=0.27]{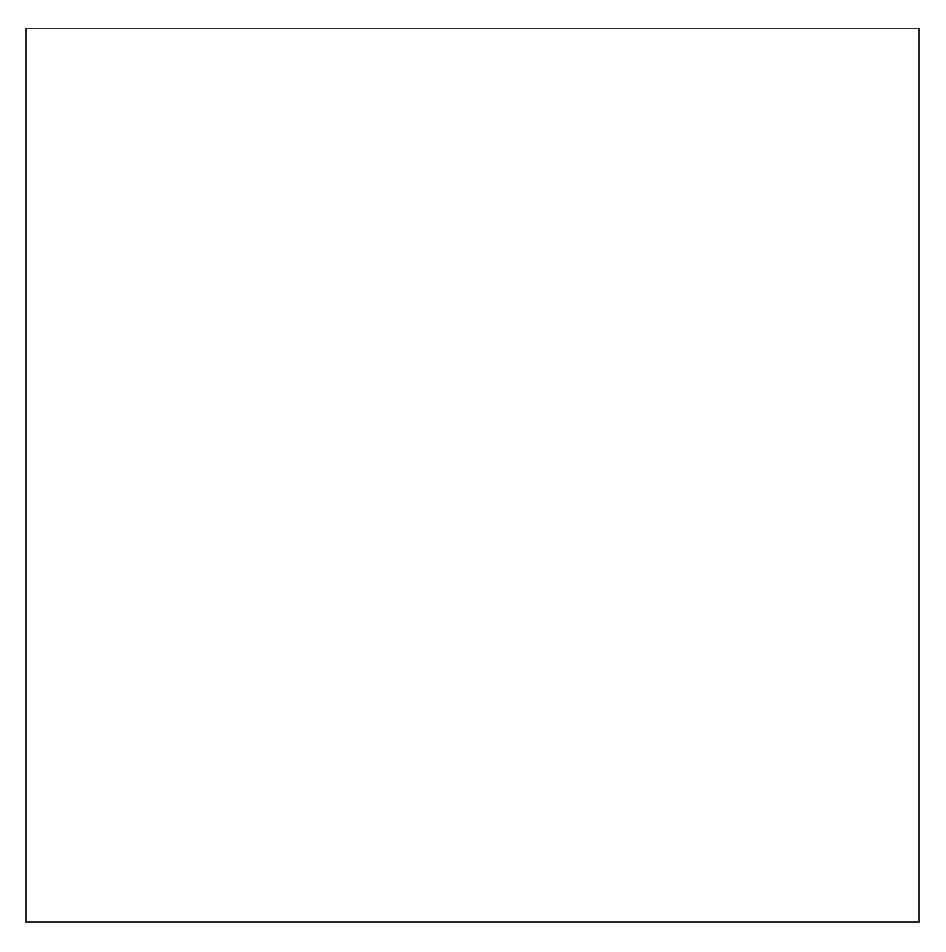}}
  \subfigure[$r=0.34, t=0$]{
    \includegraphics[scale=0.27]{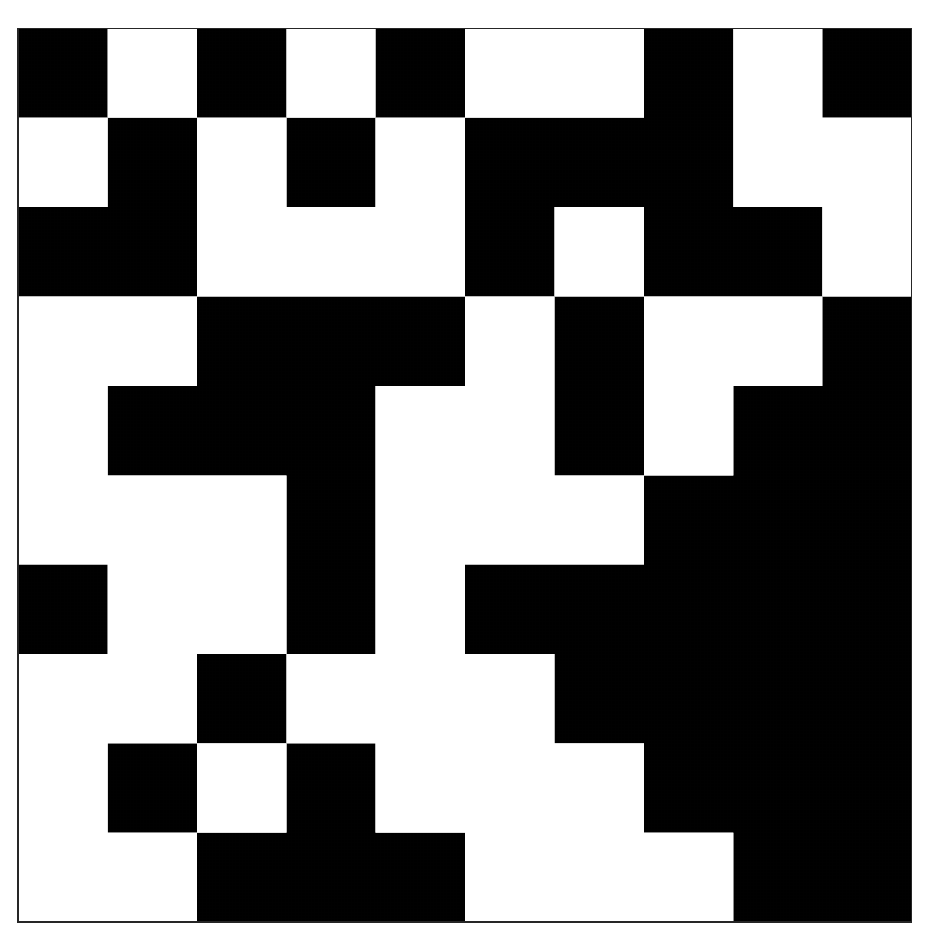}}
  \hspace{0.02in} 
  \subfigure[$r=0.34, t=5$]{
    \includegraphics[scale=0.27]{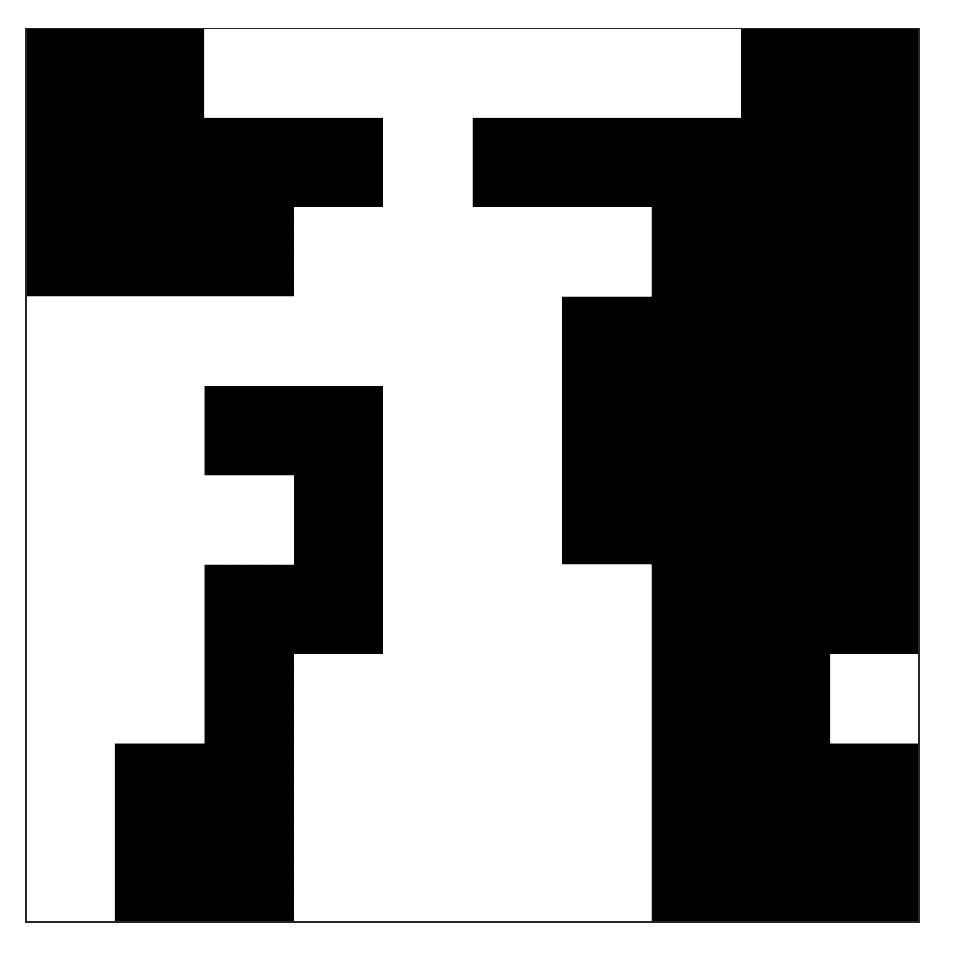}}
    \hspace{0.02in} 
  \subfigure[$b=0.34, t=18$]{
    \includegraphics[scale=0.27]{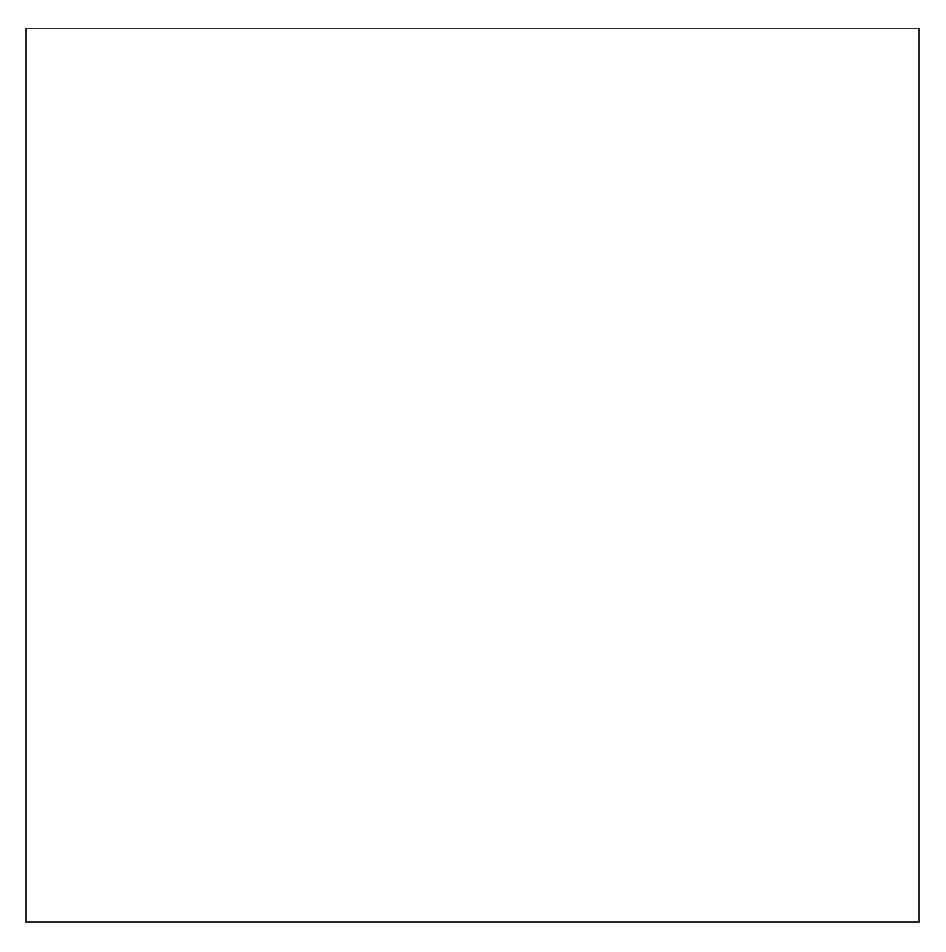}}
  \caption{The evolutionary stag hunt game with $r=0.32, 1/3, 0.34$ respectively. These systems converge fast to total cooperation.}
  \label{esh}
\end{figure}

Finally,  we simulate  evolutionary stag hunt game with the payoff matrix
\begin{table}[hbp]\label{Tab_a0}
\centering
\begin{tabular}{| c | c | c |}
\hline
      & C & D \\
\hline
      C & $1$ & $-r$  \\
\hline
      D & $r$ & $0$ \\
 \hline
\end{tabular}
\end{table}\\
satisfying $0<r<1$. By Theorem \ref{Main_result3} if $r<1/3$ and there exist four cooperation nodes forming a square at the initial time, then  evolutionary stag hunt game converges to total cooperation.
Simulations show that the evolutionary stag hunt game converges fast to total cooperation when $r=0.32, 1/3, 0.33$ (see Fig. \ref{esh}), so $r=1/3$ should not be a critical value for the convergence of evolutionary stag hunt game.

\renewcommand{\thesection}{\Roman{section}}
\section{Concluding Remarks}\label{Conclusions}
\renewcommand{\thesection}{\arabic{section}}

Since the network reciprocity has been considered one of main mechanisms accounting for the evolution of cooperative behavior, which is treated as one of the most important research topics in the 21st century,
spatial (or networked) evolutionary games are attracting more and more attention. However, due to the nature of nonlinearity, spatial evolutionary games are generally hard to analyze.
By the mehtod of  ``transforming the analysis of a stochastic system into the
design of control algorithms" first proposed by \cite{GC:17b}, this paper gives the convergence analysis of evolutionary prisoner's dilemma, evolutionary snowdrift game, and
evolutionary  stag hunt game on the two-dimensional grid.
Simulations show that our results may almost reach the critical convergence condition for the evolutionary snowdrift (or hawk-dove, chicken) game.
Also, this paper tries to solve MACC problem of  evolutionary games on the toroidal grid, and shows that for some evolutionary games (like the evolutionary prisoner's dilemma),  one fixed defection node can drive all nodes
almost surely converging to defection, while at least four fixed cooperation nodes are required to lead all nodes almost surely converging to cooperation.

There are still some problems that have not been solved. First, the sufficient and necessary condition for the convergence of evolutionary games on the grid is unknown.
 Also, we have not considered the punishment mechanism for defection.
Whether can the  punishment mechanism for defection promote cooperation? Moveover, can our results be extended to
other networks such as Erd\"{o}s-R\'{e}nyi graphs or scale-free networks? The answers to these questions may be left for us to explore in the future.

\appendices

\section{Proof of Theorem~\ref{Main_result2}}\label{sec_proof_mr2}
As same as the proof of Theorem~\ref{Main_result}, we set $\Omega^* \subseteq \{C,D\}^{N\times M}$ to be the set of strategy matrices satisfying that two adjacent nodes with different strategies have a same payoff. Similar to the proof of Theorem~\ref{Main_result}, we only need
 design a control algorithm such that $\Omega^*$ is finite-time
  reachable from any given initial $S(0)\notin\Omega^*$ under the CEG.
The algorithm is established by the following four steps:

\textbf{Step 1}: For $t\geq 0$, we repeatedly carry out the CEG with the control input (\ref{mresult_1}) until the strategies of all nodes keep unchanged.  We record the stop time as $T_1$.
For $0\leq t\leq T_1-1$,  by the relations $p_3>p_1>p_2>p_4$ and $4p_2<p_3+3 p_4$  we similarly obtain (\ref{mresult_1_2}), (\ref{mresult_1_3}) and (\ref{mresult_1_3b}), and also get that there is no isolated cooperation node or
isolated defection node at time $T_1$.
 Let $n_C(t)$ be the number of cooperation nodes at time $t$. We can similarly obtain (\ref{mresult_2}).
If $n_C(T_1)=0$, then our result is obtained. Otherwise, we carry out the following Step 2.

\textbf{Step 2}: If there exist some cooperation nodes which have exact one cooperation neighbor at time $T_1$, and have different payoff from defection neighbors,
 without loss of generality we assume
\begin{equation}\label{mr2_1}
\left\{
\begin{aligned}
&S_{i,j}(T_1)=S_{i+1,j}(T_1)=C,\\
& S_{i,j+1}(T_1)=S_{i,j-1}(T_1)=S_{i-1,j}(T_1)=D, \\
&P_{i,j}(T_1)\neq  P_{i,j+1}(T_1).
\end{aligned}\right.
\end{equation}
By (\ref{payfun}), (\ref{mr2_1}) and (\ref{mresult_1_3b}) we have
\begin{equation}\label{mr2_3}
P_{i,j}(T_1)=p_1+3p_2 > P_{i,j+1}(T_1).
\end{equation}
By the relations $p_3>p_1>p_2>p_4$ and  $p_1+p_2< p_3+p_4$ we get
$2p_3+2p_4>p_1+3p_2$, then by (\ref{payfun}) and (\ref{mr2_3}) we have
the nodes $(i,j+1),(i,j-1),(i-1,j)$ have only one cooperation neighbor at time $T_1$. Thus,
\begin{equation}\label{mr2_4}
\left\{
\begin{aligned}
& S_{i+1,j+1}(T_1)=S_{i+1,j-1}(T_1)=D, \\
& P_{i+1,j}(T_1)\leq 2p_1+2p_2.
\end{aligned}\right.
\end{equation}
Because $S_{i+1,j}(T_1)=C$ and $S_{i+1,j+1}(T_1)=D$, by (\ref{mresult_1_3b}) we have
$P_{i+1,j}(T_1)\geq P_{i+1,j+1}(T_1)$.
We will continue our discussion by the following two cases.

Case I: $P_{i+1,j+1}(T_1)=P_{i+1,j}(T_1)$. Choose $\mathcal{C}_{i,j+1}(T_1)=(i,j)$ for node $(i,j+1)$, while for any other node we choose a neighbor which has the same strategy with it as the control input of $T_1$ (We note that at time $T_1$ there is no isolated cooperation or detection node by the discussion of Step 1).  Then, by the CEG, at time $T_1+1$ the strategy of the node $(i,j+1)$ changes from $D$ to $C$, while the other nodes keep unvaried. Thus,
\begin{equation}\label{mr2_5}
P_{i+1,j+1}(T_1+1)>P_{i+1,j+1}(T_1)=P_{i+1,j}(T_1)=P_{i+1,j}(T_1).
\end{equation}
Also, by the relations $p_3>p_1>p_2>p_4$ and  $p_1+p_2< p_3+p_4$  we get
\begin{equation}\label{mr2_6}
P_{i+1,j+1}(T_1+1)\geq 2p_3+2p_4>p_1+3p_2=P_{i,j+1}(T_1).
\end{equation}
We adopt the control input (\ref{mresult_1}) at time $T_1+1$. By (\ref{mr2_5}), (\ref{mr2_6})
and (\ref{mresult_1_3}) we can get
\begin{equation}\label{mr2_7}
S_{i+1,j+1}(T_1+2)=S_{i,j+1}(T_1+2)=S_{i+1,j}(T_1+2)=D.
\end{equation}
The evolution of above process is same as Fig. \ref{LFig3}.
Combining (\ref{mr2_7}) with (\ref{mresult_1_2}) and (\ref{mr2_1}) yields
\begin{equation}\label{mr2_8}
n_{C}(T_1)> n_C(T_1+2).
\end{equation}

Case II: $P_{i+1,j+1}(T_1)<P_{i+1,j}(T_1)$. By (\ref{mr2_4}) and the relation $p_1+p_2< p_3+p_4$,
the defection node $(i+1,j+1)$ must have one cooperation neighbor nodes at time $T_1$.
We choose $\mathcal{C}_{i+1,j+1}(T_1)=(i+1,j)$ for node $(i+1,j+1)$, while for any other node we choose a neighbor which has the same strategy with it as the control input.  Then, under the CEG, at time $T_1+1$ the strategy of the node $(i+1,j+1)$ changes from $D$ to $C$, while the other nodes keep strategy unaltered. At time $T_1+1$, the cooperation node $(i+1,j+1)$ still has one cooperation neighbor, while the defection node $(i,j+1)$ has two cooperation neighbors. Thus,
 by (\ref{mr2_1}), (\ref{mr2_3}), and (\ref{payfun}),  and the relations $p_3>p_1>p_2>p_4$ and $p_1+p_2< p_3+p_4$  we have
\begin{multline}\label{mr2_9}
P_{i,j+1}(T_1+1)=2p_3+2p_4\\
>p_1+3 p_2=P_{i,j}(T_1+1)=P_{i+1,j+1}(T_1+1).
\end{multline}
We adopt the control input (\ref{mresult_1}) at time $T_1+1$. According to (\ref{mr2_9})
and (\ref{mresult_1_3}) we can get
\begin{equation*}\label{mr2_10}
S_{i,j}(T_1+2)=S_{i,j+1}(T_1+2)=S_{i+1,j+1}(T_1+2)=D,
\end{equation*}
and then (\ref{mr2_8}) still holds.
The evolution of above process is same as Fig. \ref{LFig2}.
Combining (\ref{mr2_7}) with (\ref{mresult_1_2}) and (\ref{mr2_1}) yields (\ref{mr2_8}).

For $t\geq T_1+2$  we repeatedly carry out Step 1 and the above process until  the strategies of all nodes keep unchanged under the control input (\ref{mresult_1}), and
 there is no cooperation node
which has one cooperation neighbor and has a payoff different from defection neighbors.
We record the stop time as $T_2$.
If $S(T_2)\in\Omega^*$, our result is obtained. Otherwise, we continue the following Step 3.

\textbf{Step 3}:
If there exist some cooperation nodes which have exact $2$ cooperation neighbors at time $T_2$, and have different payoff from defection neighbors,
 without loss of generality we assume $(i,j)$ has exact $2$ cooperation neighbors, and
\begin{equation}\label{mr2_11}
\left\{
\begin{aligned}
&S_{i,j}(T_2)=S_{i+1,j}(T_2)=C, S_{i,j+1}(T_2)=D, \\
&P_{i,j}(T_2)\neq  P_{i,j+1}(T_2).
\end{aligned}\right.
\end{equation}
By (\ref{payfun}), (\ref{mr2_11}) and (\ref{mresult_1_3b}) we have
\begin{equation}\label{mr2_12}
P_{i,j}(T_2)=2p_1+2p_2 > P_{i,j+1}(T_2).
\end{equation}
From (\ref{mr2_12}), and by the relations $p_3>p_1>p_2>p_4$ and $p_1+p_2< p_3+p_4$,
we can get the defection node $(i,j+1)$ has only one cooperation node at time $T_2$.
Thus,  $S_{i+1,j+1}(T_2)=D$. By (\ref{mresult_1_3b}) we have
$P_{i+1,j}(T_2)\geq P_{i+1,j+1}(T_2)$.

If $P_{i+1,j+1}(T_2)=P_{i+1,j}(T_2)$, with the similar discussion to the Case I in Step 2 we can get
\begin{equation}\label{mr2_13}
n_{C}(T_2)> n_C(T_2+2).
\end{equation}

If $P_{i+1,j+1}(T_2)<P_{i+1,j}(T_2)$, because $P_{i+1,j}(T_2)\leq 3p_1+p_2<3p_3+p_4$,
the defection node $(i+1,j+1)$ must have one or two cooperation neighbor nodes at time $T_2$.
We choose $\mathcal{C}_{i+1,j+1}(T_2)=(i+1,j)$ for node $(i+1,j+1)$, while for any other node we
 choose a neighbor which has the same strategy with it as the control input.
Then, under the CEG, at time $T_2+1$ the strategy of the node
$(i+1,j+1)$ changes from $D$ to $C$, while the other nodes keep strategy unaltered. At time $T_2+1$,
the cooperation node $(i+1,j+1)$ have at most two cooperation neighbors, while the defection node $(i,j+1)$
has two cooperation neighbors. Thus,
 by (\ref{payfun}) and the relations $p_3>p_1>p_2>p_4$  and $p_1+p_2< p_3+p_4$ we have
\begin{multline}\label{mr2_14}
P_{i,j+1}(T_2+1)=2p_3+2p_4>2p_1+2p_2\\
\geq \max\{P_{i,j}(T_2+1),P_{i+1,j+1}(T_2+1)\}.
\end{multline}
We adopt the control input (\ref{mresult_1}) at time $T_2+1$. According to (\ref{mr2_14})
and (\ref{mresult_1_3}) we can get
\begin{equation*}\label{mr2_15}
S_{i,j}(T_2+2)=S_{i,j+1}(T_2+2)=S_{i+1,j+1}(T_2+2)=D,
\end{equation*}
and then (\ref{mr2_13}) still holds.
Fig. \ref{LFig2} shows the evolution of nodes' strategies for the case when $P_{i+1,j+1}(T_2)<P_{i+1,j}(T_2)$.

For $t\geq T_2+2$  we repeatedly carry out Steps 1-2 and the above process until  the strategies of all nodes keep unchanged under the control input (\ref{mresult_1}),
and there is no cooperation node
which has one or two cooperation neighbors and has a payoff different from defection neighbors.
We record the stop time as $T_3$.
If $S(T_3)\in\Omega^*$, our result is obtained. Otherwise, we continue the following Step 4.

\textbf{Step 4}: By the similar method to the Step 3 in the proof of Theorem \ref{Main_result} we can obtain
\begin{equation}\label{mr2_16}
n_{C}(T_3)> n_C(T_3+2).
\end{equation}
Here we note that with the relations $p_3>p_1>p_2>p_4$ and $p_1+p_2<p_3+p_4$, the inequations
$2 p_3+2p_4>2 p_1+2p_2$ in (\ref{mrb_10_2}), $p_3-p_4>0$ in (\ref{mrb_10_3}), and $3p_3+p_4>3p_1+p_2$
in (\ref{mresult_2_1b}) still holds.

For $t\geq T_3+2$  we repeatedly carry out Steps 1-3 and the above process until $\Omega^*$ is reached. Let $T_4$ be the stop time. By (\ref{mresult_2}), (\ref{mr2_13}), (\ref{mr2_8}) and (\ref{mr2_16}) we have $T_4\leq 2(NM-1)$.

\ifCLASSOPTIONcaptionsoff
  \newpage
\fi

\bibliographystyle{IEEEtran}
\bibliography{alias,FB,Main}

\end{document}